\documentclass[preprint,12pt]{elsarticle}

\usepackage{amsfonts}
\usepackage{graphicx}
\usepackage{epsfig}
\usepackage{geometry}
\usepackage{tikz}
\usetikzlibrary{trees,positioning,shapes}
\usepackage{subcaption}
\usepackage{amsmath}
\usepackage{amsthm}

\usepackage{xcolor, colortbl}

\newcommand{\inblue}[1]{\textcolor{blue}{#1}}

\newcommand{\inred}[1]{\textcolor{red}{#1}}
\newcommand{\ingray}[1]{\textcolor{gray}{#1}}
\newcommand{\revision}[1]{\textcolor{black}{#1}}
\newcommand{\inbrown}[1]{\textcolor{brown}{#1}}

\newcommand{\inmagenta}[1]{\textcolor{magenta}{#1}}
\newcommand{\incyan}[1]{\textcolor{cyan}{#1}}

\newcommand{\Mat}[1]{\mathbf{#1}}

\newcommand{\B}{\mathcal{B}}

\newcommand{\M}{\Mat{M}}

\newcommand{\grad}{\nabla}

\DeclareMathOperator*{\Div}{div}

\begin{document}

\begin{frontmatter}

\title{Alternating directions implicit higher-order \\ finite element method for simulations \\ of time-dependent electromagnetic wave propagation \\ in non-regular biological tissues}

\author{Maciej Paszy\'nski$^{1}$, Marcin \L{}o\'s$^{1}$, Judit Mu\~noz-Matute$^{2,3}$}

\address{$^{(1)}$Institute of Computer Science, AGH University of Science and Technology, Krakow, Poland \\
$^{(2)}$Basque Center for Applied Mathematics (BCAM), Bilbao, Spain\\
$^{(3)}$Oden Institute for Computational Engineering and Sciences, The University of Texas at Austin, USA\\}

\begin{abstract}
We focus on non-stationary Maxwell equations defined on a regular patch of elements as considered in the isogeometric analysis (IGA). We apply the time-integration scheme following the ideas developed by the finite difference community \cite{maxwell1} to derive a weak formulation resulting in discretization with Kronecker product matrices. 
We take the tensor product structure of the computational patch of elements from the IGA framework as an advantage, allowing for linear computational cost factorization in every time step.
We design our solver to target simulations of electromagnetic waves propagations in non-regular biological tissues.  We show that the linear cost of the alternating direction solver 
is preserved when we arbitrarily vary material data coefficients across the computational domain. We verify the solver using the manufactured solution and the problem of propagation of electromagnetic waves on the human head.
\end{abstract}
	
%\begin{keyword}
%isogeometric analysis \sep alternating directions \sep residual minimization \sep time-dependent Maxwell problem \sep Kronecker product \sep linear computational cost solver
%\end{keyword}

\end{frontmatter}

\section*{Introduction}

In this paper, we introduce a fast solver for non-stationary simulations of propagation of  electromagnetic waves over non-regular biological tissues,  with the following unique combination of features:
\begin{enumerate}
\item  Linear computational cost $O(N)$ of the direct solution. 
%even for non-regular material data.
\item Unconditional stability of the implicit time integration scheme. 
\item Second order accurate time integration scheme.
\end{enumerate}
%We believe this new paradigm may have an impact on the way the computational mechanics community runs simulations. 
It mixes benefits of the state-of-the-art modern methods, the Isogeometric Finite Element Method (IGA-FEM) \cite{p1}, and Alternating Direction Implicit solvers (ADI) \cite{Sportisse,maxwell1}.
IGA-FEM utilizes higher-order and continuity basis functions to obtain smooth and continuous approximation of the solution vector fields.
Splitting methods modify the original linear systems of equations seeking to reduce computation costs. For instance, operator splitting methods decrease the dimension of the matrices and handle implicit time marching efficiently \cite{Sportisse,ADS1,ADS2}.
We exploit the tensor product structure of the discretization to represent the system matrices as Kronecker products of
one-dimensional matrices. This reinterpretation of the algebraic system allows us to design simple
approximations that deliver linear computational cost. These approximations are based on the alternating direction method. The results detailed in \cite{Gao1,Gao2,Los1,Los2,iGRM1,iGRM2} prove that alternating direction splitting solvers based on tensor-products result in linear
computational cost for every time step.
It is a common misunderstanding that the direction splitting solvers are limited to simple geometries. On the contrary, they can be applied to discretizations in extremely complicated geometries, as described in \cite{Minev}. However, this requires the development of problem-specific methods and implementations.
In particular, we show that our solver can be applied when we arbitrarily vary material data coefficients across the computational domain.

In this paper, we focus on the non-stationary Maxwell problem.
We employ the B-spline basis functions from isogeometric analysis (IGA) for the discretization in space.

The direction splitting method was applied to the Maxwell problem in the context of the finite difference method ~\cite{ maxwell1,maxwell2}. 
This paper employs the finite element method with higher continuity B-spline basis functions for discretization.

We design our solver for electromagnetic waves propagations in non-regular biological tissues, and we utilize the MRI scan of the human head to illustrate the concept.
We show that the alternating direction splitting algorithm can be applied when we arbitrarily vary material data coefficients across the computational domain, including the tissue, skull, and air. We verify our solver using the manufactured solution technique. Finally, we summarize the paper with a numerical example of the propagation of electromagnetic waves over the human head.

The structure of the paper is the following. We start with the introduction of the direction splitting method for the non-stationary Maxwell problem, following \cite{maxwell1,maxwell2}. Next, we introduce the variational formulations with B-spline basis functions on tensor product grids preserving the Kronecker product structure of matrices. Later we verify the method using a numerical example with the manufactured solution. Finally, we introduce the algorithm for incorporating non-regular biological tissues of the human head.
We summarize the paper with the numerical experiment of the propagation of electromagnetic waves over the human head. Our algorithm is summarized in the appendix, and the main results are summarized in the conclusions section.

\section*{Alternating directions splitting for Maxwell equations}

Let us consider the time-depedent Maxwell equations on domain~$\Omega=(0,1)^3$:
\begin{equation}
\begin{aligned}
\frac{\partial {\bf E}}{\partial t}(t) &= \frac{1}{\epsilon} \grad \times {\bf H}(t),
&
\frac{\partial {\bf H}}{\partial t}(t)&=-\frac{1}{\mu} \grad \times {\bf E}(t)
&& (x, t) \in \Omega \times (0, T) \\
\Div \epsilon {\bf E}(t) &= 0,
&
\Div{\mu} {\bf H}(t) &=0
&& (x, t) \in \Omega \times (0, T) \\
{\bf E}(t)  \times {\bf n} &=0,
&
{\bf H}(t)  \cdot {\bf n} &=0
&&
(x, t) \in \partial\Omega \times (0, T) \\
{\bf E}(x,0)&={\bf E_0}(x),
&
{\bf H}(x,0)&={\bf H_0}(x)
&&
x \in \Omega
\end{aligned}
\end{equation}
where ${\bf E}(x,t)$ is the electric field and ${\bf H}(x,t)$ is the magnetic field. Here, ${\bf E_0}\in L^2(\Omega)^3$ and ${\bf H_0}\in L^2(\Omega)^3$ are initial states. The permittivity $\epsilon\in L^{\infty}(\Omega)$ and the permeability $\mu\in L^{\infty}(\Omega)$ are given functions assumed to be constant in time, and they fullfil $\epsilon(x)\geq \delta >0$, and $\mu(x) \geq \delta > 0$.

We employ the same time-integration scheme as in \cite{maxwell1}. For that, we first split the curl operator as
\begin{equation}
\grad \times = 
\begin{bmatrix}0 &  -\frac{\partial}{\partial x_3} & \frac{\partial}{\partial x_2} \\ 
\frac{\partial}{\partial x_3} & 0 & -\frac{\partial}{\partial x_1}\\  -  \frac{\partial}{\partial x_2}  & \frac{\partial}{\partial x_1} & 0\\ 
\end{bmatrix} = \begin{bmatrix}0 & 0 & \frac{\partial}{\partial x_2} \\ 
\frac{\partial}{\partial x_3} & 0 & 0\\  0 & \frac{\partial}{\partial x_1} & 0\\ 
\end{bmatrix}-
\begin{bmatrix}0 & \frac{\partial}{\partial x_3} & 0  \\ 
0 & 0& \frac{\partial}{\partial x_1} \\  \frac{\partial}{\partial x_2} & 0 & 0\\ 
\end{bmatrix}=C_1-C_2
\end{equation}
and we consider the following time-marching scheme that consists in two substeps
\begin{equation}
{\bf E}^{n+\frac{1}{2}}={\bf E}^n - \frac{\tau }{2 \epsilon} C_2 {\bf H}^n + \frac{\tau }{2 \epsilon} C_1  {\bf H}^{n+\frac{1}{2}}, \qquad
{\bf H}^{n+\frac{1}{2}}={\bf H}^n - \frac{\tau }{2 \mu} C_1 {\bf E}^n + \frac{\tau }{2 \mu} C_2  {\bf E}^{n+\frac{1}{2}}
\end{equation}
\begin{equation}
{\bf E}^{n+1}={\bf E}^{n+\frac{1}{2}} + \frac{\tau }{2 \epsilon} C_1 {\bf H}^{n+\frac{1}{2}} - \frac{\tau }{2 \epsilon} C_2  {\bf H}^{n+1}, \qquad
{\bf H}^{n+1}={\bf H}^{n+\frac{1}{2}} + \frac{\tau }{2 \mu} C_2 {\bf E}^{n+\frac{1}{2}} - \frac{\tau }{2 \mu} C_1  {\bf E}^{n+1}
\end{equation}

Substituting the second equation in the first one in both substeps lead to
\begin{equation}
\left(I-  \frac{\tau^2 }{4 \epsilon} C_1 \mu^{-1} C_2\right) {\bf E}^{n+\frac{1}{2}}=
{\bf E}^n + \frac{\tau }{2 \epsilon} (C_1-C_2) {\bf H}^n - \frac{\tau^2 }{4\epsilon} C_1\mu^{-1}C_1 {\bf E}^n\label{scheme1a}
\end{equation} 
\begin{equation}
{\bf H}^{n+\frac{1}{2}}={\bf H}^n - \frac{\tau }{2 \mu} C_1 {\bf E}^n + \frac{\tau }{2 \mu} C_2  {\bf E}^{n+\frac{1}{2}}\label{scheme1b}
\end{equation}
and
\begin{equation}
\left(I-  \frac{\tau^2 }{4 \epsilon} C_2 \mu^{-1} C_1\right) {\bf E}^{n+1}=
{\bf E}^{n+\frac{1}{2}} + \frac{\tau }{2 \epsilon} (C_1-C_2) {\bf H}^{n+\frac{1}{2}} - \frac{\tau^2 }{4\epsilon} C_2\mu^{-1}C_2 {\bf E}^{n+\frac{1}{2}}\label{scheme1c}
\end{equation}
\begin{equation}
{\bf H}^{n+1}={\bf H}^{n+\frac{1}{2}} \revision{+} \frac{\tau }{2 \mu} C_2 {\bf E}^{n+\frac{1}{2}} \revision{-} \frac{\tau }{2 \mu} C_1  {\bf E}^{n+1}\label{scheme1d}
\end{equation}

Since
\begin{equation}
\begin{split}
C_1\mu^{-1}C_2&= \begin{bmatrix} \frac{\partial}{\partial x_2}\mu^{-1}\frac{\partial}{\partial x_2}& 0 & 0 \\ 
0 &\frac{\partial}{\partial x_3} \mu^{-1} \frac{\partial}{\partial x_3}  & 0\\  0 & 0 & \frac{\partial}{\partial x_1}\mu^{-1}\frac{\partial}{\partial x_1} \\ 
\end{bmatrix}\\
C_2\mu^{-1}C_1&= \begin{bmatrix} \frac{\partial}{\partial x_3}\mu^{-1}\frac{\partial}{\partial x_3}& 0 & 0 \\ 
0 &\frac{\partial}{\partial x_1} \mu^{-1} \frac{\partial}{\partial x_1}  & 0\\  0 & 0 & \frac{\partial}{\partial x_2}\mu^{-1}\frac{\partial}{\partial x_2} \\ 
\end{bmatrix}
\end{split}
\end{equation}
as well as
\begin{equation}
\begin{split}
C_1\mu^{-1}C_1&= \begin{bmatrix} 0 & \frac{\partial}{\partial x_2}\mu^{-1}\frac{\partial}{\partial x_1}& 0 \\ 
0 & 0 & \frac{\partial}{\partial x_3} \mu^{-1} \frac{\partial}{\partial x_2}  \\  \frac{\partial}{\partial x_1}\mu^{-1}\frac{\partial}{\partial x_3} & 0 & 0 \\ 
\end{bmatrix}\\
C_2\mu^{-1}C_2&= \begin{bmatrix} 0 & 0 & \frac{\partial}{\partial x_3}\mu^{-1}\frac{\partial}{\partial x_1} \\ 
\frac{\partial}{\partial x_1} \mu^{-1} \frac{\partial}{\partial x_2} & 0 & 0 \\ 0 &  \frac{\partial}{\partial x_2}\mu^{-1}\frac{\partial}{\partial x_3} & 0 \\ 
\end{bmatrix}
\end{split}
\end{equation}
we obtain
\begin{equation}\label{strong1}
\begin{split}
{\bf E}^{n+\frac{1}{2}}-
 \frac{\tau^2 }{4 \epsilon}\begin{bmatrix}\frac{\partial}{\partial x_2}\mu^{-1}\frac{\partial}{\partial x_2}& 0 & 0 \\ 
0 &\frac{\partial}{\partial x_3} \mu^{-1} \frac{\partial}{\partial x_3}  & 0\\  0 & 0 &  \frac{\partial}{\partial x_1}\mu^{-1}\frac{\partial}{\partial x_1} \\ 
\end{bmatrix} {\bf E}^{n+\frac{1}{2}}= \\
{\bf E}^n + \frac{\tau }{2 \epsilon} \begin{bmatrix}0 &  -\frac{\partial}{\partial x_3} & \frac{\partial}{\partial x_2} \\ 
\frac{\partial}{\partial x_3} & 0 & -\frac{\partial}{\partial x_1}\\  -  \frac{\partial}{\partial x_2}  & \frac{\partial}{\partial x_1} & 0\\ 
\end{bmatrix}  {\bf H}^n \\ - \frac{\tau^2 }{4\epsilon} \begin{bmatrix} 0 & \frac{\partial}{\partial x_2}\mu^{-1}\frac{\partial}{\partial x_1}& 0 \\ 
0 & 0 & \frac{\partial}{\partial x_3} \mu^{-1} \frac{\partial}{\partial x_2}  \\  \frac{\partial}{\partial x_1}\mu^{-1}\frac{\partial}{\partial x_3} & 0 & 0 
\end{bmatrix} {\bf E}^n
\end{split}
\end{equation}
\begin{equation}\label{strong2}
{\bf H}^{n+\frac{1}{2}}={\bf H}^n - \frac{\tau }{2 \mu} \begin{bmatrix}0 & 0 & \frac{\partial}{\partial x_2} \\ 
\frac{\partial}{\partial x_3} & 0 & 0\\  0 & \frac{\partial}{\partial x_1} & 0
\end{bmatrix} {\bf E}^n +\frac{\tau }{2 \mu} \begin{bmatrix}0 & \frac{\partial}{\partial x_3} & 0  \\ 
0 & 0& \frac{\partial}{\partial x_1} \\  \frac{\partial}{\partial x_2} & 0 & 0 
\end{bmatrix} {\bf E}^{n+\frac{1}{2}}
\end{equation}
and
\begin{equation}\label{strong3}
\begin{split}
{\bf E}^{n+1}-
\frac{\tau^2 }{4 \epsilon}\begin{bmatrix}   \frac{\partial}{\partial x_3}\mu^{-1}\frac{\partial}{\partial x_3}& 0 & 0 \\ 
0 & \frac{\partial}{\partial x_1} \mu^{-1} \frac{\partial}{\partial x_1}  & 0\\  0 & 0 & \frac{\partial}{\partial x_2}\mu^{-1}\frac{\partial}{\partial x_2} \\ 
\end{bmatrix} {\bf E}^{n+1}=\\
{\bf E}^{n+\frac{1}{2}}+ \frac{\tau }{2 \epsilon} \begin{bmatrix}0 &  -\frac{\partial}{\partial x_3} & \frac{\partial}{\partial x_2} \\ 
\frac{\partial}{\partial x_3} & 0 & -\frac{\partial}{\partial x_1}\\  -  \frac{\partial}{\partial x_2}  & \frac{\partial}{\partial x_1} & 0  
\end{bmatrix}  {\bf H}^{n+\frac{1}{2}} \\ - \frac{\tau^2 }{4\epsilon} \begin{bmatrix} 0 & 0 & \frac{\partial}{\partial x_3}\mu^{-1}\frac{\partial}{\partial x_1} \\ 
\frac{\partial}{\partial x_1} \mu^{-1} \frac{\partial}{\partial x_2} & 0 & 0 \\ 0 &  \frac{\partial}{\partial x_2}\mu^{-1}\frac{\partial}{\partial x_3} & 0 \\ 
\end{bmatrix}{\bf E}^{n+\frac{1}{2}}
\end{split}
\end{equation}
\begin{equation}\label{strong4}
{\bf H}^{n+1}={\bf H}^{n+\frac{1}{2}} \revision{+} \frac{\tau }{2 \mu} \begin{bmatrix}0 & \frac{\partial}{\partial x_3} & 0  \\ 
0 & 0& \frac{\partial}{\partial x_1} \\  \frac{\partial}{\partial x_2} & 0 & 0
\end{bmatrix} {\bf E}^{n+\frac{1}{2}} - \frac{\tau }{2 \mu} \begin{bmatrix}0 & 0 & \frac{\partial}{\partial x_2} \\ 
\frac{\partial}{\partial x_3} & 0 & 0\\  0 & \frac{\partial}{\partial x_1} & 0
\end{bmatrix}  {\bf E}^{n+1} 
\end{equation}

\section*{Variational formulation}
{In this section, we introduce a variational formulation of equations (\ref{strong1})-(\ref{strong4}).We denote by $(\cdot,\cdot)$ both the usual inner products in $L^2(\Omega)$ and $L^2(\Omega)^3$, i.,e.
$$(\bold H,\bold V)=(H_1,V_1)+(H_2,V_2)+(H_3,V_3)$$
We consider for the moment that $\mu$ and $\epsilon$ are constant. We multiply the equations by suitable test functions $\bf V$, integrate in space and integrate by parts the second order terms
\begin{equation}\label{weak1}
\begin{split}
(E_1^{n+\frac{1}{2}},V_1)+(E_2^{n+\frac{1}{2}},V_2)+(E_3^{n+\frac{1}{2}},V_3)\\
+\frac{\tau^2}{4\epsilon\mu}\left[\left(\frac{\partial}{\partial x_2}E_1^{n+\frac{1}{2}},\frac{\partial}{\partial x_2}V_1\right)+\left(\frac{\partial}{\partial x_3}E_2^{n+\frac{1}{2}},\frac{\partial}{\partial x_3}V_2\right)+\left(\frac{\partial}{\partial x_1}E_3^{n+\frac{1}{2}},\frac{\partial}{\partial x_1}V_3\right)\right]\\
=(E_1^{n},V_1)+(E_2^{n},V_2)+(E_3^{n},V_3)\\
+\frac{\tau}{2\epsilon}\left[\left(-\frac{\partial}{\partial x_3}H_2^n+\frac{\partial}{\partial x_2}H_3^n,V_1\right)+\left(\frac{\partial}{\partial x_3}H_1^n-\frac{\partial}{\partial x_1}H_3^n,V_2\right)\right.\\
\left.+\left(-\frac{\partial}{\partial x_2}H_1^n+\frac{\partial}{\partial x_1}H_2^n,V_3\right)\right]\\
+\frac{\tau^2}{4\epsilon\mu}\left[\left(\frac{\partial}{\partial x_1}E_2^{n},\frac{\partial}{\partial x_2}V_1\right)+\left(\frac{\partial}{\partial x_2}E_3^{n},\frac{\partial}{\partial x_3}V_2\right)+\left(\frac{\partial}{\partial x_3}E_1^{n},\frac{\partial}{\partial x_1}V_3\right)\right]\\
\end{split}
\end{equation}
\vspace{0.3cm}
\begin{equation}\label{weak2}
\begin{split}
(H_1^{n+\frac{1}{2}},V_1)+(H_2^{n+\frac{1}{2}},V_2)+(H_3^{n+\frac{1}{2}},V_3)=(H_1^{n},V_1)+(H_2^{n},V_2)+(H_3^{n},V_3)\\
-\frac{\tau}{2\mu}\left[\left(\frac{\partial}{\partial x_2}E_3^{n},V_1\right)+\left(\frac{\partial}{\partial x_3}E_1^{n},V_2\right)+\left(\frac{\partial}{\partial x_1}E_2^{n},V_3\right)\right]\\
+\frac{\tau}{2\mu}\left[\left(\frac{\partial}{\partial x_3}E_2^{n+\frac{1}{2}},V_1\right)+\left(\frac{\partial}{\partial x_1}E_3^{n+\frac{1}{2}},V_2\right)+\left(\frac{\partial}{\partial x_2}E_1^{n+\frac{1}{2}},V_3\right)\right]\\
\end{split}
\end{equation}
\vspace{0.3cm}
\begin{equation}\label{weak3}
\begin{split}
(E_1^{n+1},V_1)+(E_2^{n+1},V_2)+(E_3^{n+1},V_3)\\
+\frac{\tau^2}{4\epsilon\mu}\left[\left(\frac{\partial}{\partial x_3}E_1^{n+1},\frac{\partial}{\partial x_3}V_1\right)+\left(\frac{\partial}{\partial x_1}E_2^{n+1},\frac{\partial}{\partial x_1}V_2\right)+\left(\frac{\partial}{\partial x_2}E_3^{n+1},\frac{\partial}{\partial x_2}V_3\right)\right]\\
=(E_1^{n+\frac{1}{2}},V_1)+(E_2^{n+\frac{1}{2}},V_2)+(E_3^{n+\frac{1}{2}},V_3)\\
+\frac{\tau}{2\epsilon}\left[\left(-\frac{\partial}{\partial x_3}H_2^{n+\frac{1}{2}}+\frac{\partial}{\partial x_2}H_3^{n+\frac{1}{2}},V_1\right)+\left(\frac{\partial}{\partial x_3}H_1^{n+\frac{1}{2}}-\frac{\partial}{\partial x_1}H_3^{n+\frac{1}{2}},V_2\right)\right.\\
+\left.\left(-\frac{\partial}{\partial x_2}H_1^{n+\frac{1}{2}}+\frac{\partial}{\partial x_1}H_2^{n+\frac{1}{2}},V_3\right)\right]\\
+\frac{\tau^2}{4\epsilon\mu}\left[\left(\frac{\partial}{\partial x_1}E_3^{n+\frac{1}{2}},\frac{\partial}{\partial x_3}V_1\right)+\left(\frac{\partial}{\partial x_2}E_1^{n+\frac{1}{2}},\frac{\partial}{\partial x_1}V_2\right)+\left(\frac{\partial}{\partial x_3}E_2^{n+\frac{1}{2}},\frac{\partial}{\partial x_2}V_3\right)\right]\\
\end{split}
\end{equation}
\vspace{0.3cm}
\begin{equation}\label{weak4}
\begin{split}
(H_1^{n+1},V_1)+(H_2^{n+1},V_2)+(H_3^{n+1},V_3)=(H_1^{n+\frac{1}{2}},V_1)+(H_2^{n+\frac{1}{2}},V_2)+(H_3^{n+\frac{1}{2}},V_3)\\
+\frac{\tau}{2\mu}\left[\left(\frac{\partial}{\partial x_3}E_2^{n+\frac{1}{2}},V_1\right)+\left(\frac{\partial}{\partial x_1}E_3^{n+\frac{1}{2}},V_2\right)+\left(\frac{\partial}{\partial x_2}E_1^{n+\frac{1}{2}},V_3\right)\right]\\
-\frac{\tau}{2\mu}\left[\left(\frac{\partial}{\partial x_2}E_3^{n+1},V_1\right)+\left(\frac{\partial}{\partial x_3}E_1^{n+1},V_2\right)+\left(\frac{\partial}{\partial x_1}E_2^{n+1},V_3\right)\right]\\
\end{split}
\end{equation}
}

{Finally, after discretizing in space, we obtain the matrix form of equations (\ref{weak1})-(\ref{weak4})
\begin{equation}\label{discrete1}
\left({\bf M}+\frac{\tau^2}{4\epsilon\mu}{\bf K}_1\right){\bf E}^{n+\frac{1}{2}}={\bf M}{\bf E}^n+\frac{\tau}{2\epsilon}{\bf C}{\bf H}^n+\frac{\tau^2}{4\epsilon\mu}{\bf R}_1{\bf E}^n
\end{equation}
\begin{equation}\label{discrete2}
{\bf M}{\bf H}^{n+\frac{1}{2}}={\bf M}{\bf H}^n-\frac{\tau}{2\mu}{\bf C}_1{\bf E}^n+\frac{\tau}{2\mu}{\bf C}_2{\bf E}^{n+\frac{1}{2}}
\end{equation}
\begin{equation}\label{discrete3}
\left({\bf M}+\frac{\tau^2}{4\epsilon\mu}{\bf K}_2\right){\bf E}^{n+1}={\bf M}{\bf E}^{n+\frac{1}{2}}+\frac{\tau}{2\epsilon}{\bf C}{\bf H}^{n+\frac{1}{2}}+\frac{\tau^2}{4\epsilon\mu}{\bf R}_2{\bf E}^{n+\frac{1}{2}}
\end{equation}
\begin{equation}\label{discrete4}
{\bf M}{\bf H}^{n+1}={\bf M}{\bf H}^{n+\frac{1}{2}}+\frac{\tau}{2\mu}{\bf C}_2{\bf E}^{n+\frac{1}{2}}+\frac{\tau}{2\mu}{\bf C}_1{\bf E}^{n+1}
\end{equation}
where the matrices are defined as follows
\begin{equation}
\begin{split}
{\bf M}&=\begin{bmatrix}M_{x_1}\otimes M_{x_2}\otimes M_{x_3}&0&0\\0&M_{x_1} \otimes M_{x_2}\otimes M_{x_3}&0\\0&0&M_{x_1} \otimes M_{x_2}\otimes M_{x_3}\end{bmatrix}\\
{\bf K}_1&=\begin{bmatrix}M_{x_1}\otimes S_{x_2}\otimes M_{x_3}&0&0\\0&M_{x_1} \otimes M_{x_2}\otimes S_{x_3}&0\\0&0&S_{x_1} \otimes M_{x_2}\otimes M_{x_3}\end{bmatrix}\\
{\bf K}_2&=\begin{bmatrix}M_{x_1}\otimes M_{x_2}\otimes S_{x_3}&0&0\\0&S_{x_1} \otimes M_{x_2}\otimes M_{x_3}&0\\0&0&M_{x_1} \otimes S_{x_2}\otimes M_{x_3}\end{bmatrix}\\
{\bf R}_1&=\begin{bmatrix}0&A_{x_1}\otimes B_{x_2}\otimes M_{x_3}&0\\0&0&M_{x_1}\otimes A_{x_2}\otimes B_{x_3}\\B_{x_1}\otimes M_{x_2}\otimes A_{x_3}&0&0\end{bmatrix}\\
{\bf R}_2&=\begin{bmatrix}0&0&A_{x_1}\otimes M_{x_2}\otimes B_{x_3}\\B_{x_1}\otimes A_{x_2}\otimes M_{x_3}&0&0\\0&M_{x_1}\otimes B_{x_2}\otimes A_{x_3}&0\end{bmatrix}\\
{\bf C}_1&=\begin{bmatrix}0&0&M_{x_1}\otimes A_{x_2}\otimes M_{x_3}\\M_{x_1}\otimes M_{x_2}\otimes A_{x_3}&0&0\\0&A_{x_1}\otimes M_{x_2}\otimes M_{x_3}&0\end{bmatrix}\\
{\bf C}_2&=\begin{bmatrix}0&M_{x_1}\otimes M_{x_2}\otimes A_{x_3}&0\\0&0&A_{x_1}\otimes M_{x_2}\otimes M_{x_3}\\M_{x_1}\otimes M_{x_2}\otimes A_{x_3}&0&0\end{bmatrix}\\
{\bf C}&={\bf C}_1-{\bf C}_2\\
\end{split}
\end{equation}}

{where $M_{x_1},M_{x_2},M_{x_3}$ are 1D mass matrices, $S_{x_1},S_{x_2},S_{x_3}$ are 1D stiffness matrices, $A_{x_1},A_{x_2},A_{x_3}$ are 1D advection matrices with the derivatives in the trial functions, and $B_{x_1},B_{x_2},B_{x_3}$ are 1D advection matrices with the derivatives in the trial functions. In the method presented here, we obtain} Kronecker product matrices on the left-hand sides, which can be factorized in a linear cost in every time step of the time-dependent simulation. Thus, we have generalized the result from \cite{maxwell1} into isogeometric finite element method computations performed over the patch of elements. We still deliver linear computational cost solver with implicit time integration scheme for non-stationary Maxwell equations. Additionally, we provide higher-order and{high-continuity discretizations in space} as available with B-spline basis functions.

\section*{Numerical code verification with manufactured solution}

For $\Omega=(0,1)^3$, for $\epsilon=1$ and $\mu=1$ we define
\begin{equation}
u^1_{\kappa,\lambda}(x,t)=\begin{bmatrix} \sin (\kappa \pi x_2)  \sin (\lambda \pi x_3) \cos(\sqrt{\kappa^2+\lambda^2}\pi t) \\ 0 \\ 0 \\ 0 \\ 
- \frac{\lambda}{\sqrt{\kappa^2+\lambda^2}} \sin (\kappa \pi x_2)  \cos (\lambda \pi x_3) \sin(\sqrt{\kappa^2+\lambda^2}\pi t)  \\
 \frac{\kappa}{\sqrt{\kappa^2+\lambda^2}} \cos (\kappa \pi x_2)  \sin (\lambda \pi x_3) \sin(\sqrt{\kappa^2+\lambda^2}\pi t)  
\end{bmatrix}
\end{equation}
\begin{equation}
u^2_{\kappa,\lambda}(x,t)=\begin{bmatrix} 0 \\ \sin (\kappa \pi x_1)  \sin (\lambda \pi x_3) \cos(\sqrt{\kappa^2+\lambda^2}\pi t) \\ 0 \\ 
- \frac{\lambda}{\sqrt{\kappa^2+\lambda^2}} \sin (\kappa \pi x_1)  \cos (\lambda \pi x_3) \sin(\sqrt{\kappa^2+\lambda^2}\pi t)  \\ 0 \\
 \frac{\kappa}{\sqrt{\kappa^2+\lambda^2}} \cos (\kappa \pi x_1)  \sin (\lambda \pi x_3) \sin(\sqrt{\kappa^2+\lambda^2}\pi t)   
\end{bmatrix}
\end{equation}
\begin{equation}
u^3_{\kappa,\lambda}(x,t)=\begin{bmatrix} 0 \\ 0 \\ \sin (\kappa \pi x_1)  \sin (\lambda \pi x_2) \cos(\sqrt{\kappa^2+\lambda^2}\pi t) \\ 
- \frac{\lambda}{\sqrt{\kappa^2+\lambda^2}} \sin (\kappa \pi x_1)  \cos (\lambda \pi x_2) \sin(\sqrt{\kappa^2+\lambda^2}\pi t)  \\ 
 \frac{\kappa}{\sqrt{\kappa^2+\lambda^2}} \cos (\kappa \pi x_1)  \sin (\lambda \pi x_2) \sin(\sqrt{\kappa^2+\lambda^2}\pi t)   \\ 0
\end{bmatrix}
\end{equation}
for $\kappa,\lambda \in \mathbb{Z}, \kappa,\lambda  \neq 0$. 

The first manufactured solution function is
\begin{equation} 
{\bf u}_A(x,t)=\gamma  u^1_{1,1}(x,t) + 2\gamma u^2_{1,1}(x,t) +3\gamma u^3_{1,1}(x,t)
\end{equation} 
Notice that ${\bf u}_A$ has six components, where the first three components correspond to ${\bf E}$ and the last three components to ${\bf H}$.
The parameter $\gamma$ is selected in such a way that $ \|{\bf u}_A(x,0)\|_{L^2(\Omega)}=1$. Since
%In other words, $ \|{\bf u}_A(x,0)\|^2_{L^2(\Omega)}=1$. 
\begin{equation}
u^1_{1,1}(x,0)=\begin{bmatrix} \sin (\pi x_2)  \sin (\pi x_3)  \\ 0 \\ 0 \\ 0 \\ 
0 \\
0
\end{bmatrix} 
\end{equation}
\begin{equation}
u^2_{1,1}(x,0)=\begin{bmatrix} 0 \\ \sin (\pi x_1)  \sin ( \pi x_3)  \\ 0 \\ 
0 \\ 0 \\ 0 
\end{bmatrix} 
\end{equation}
\begin{equation}
u^3_{1,1}(x,0)=\begin{bmatrix} 0 \\ 0 \\ \sin ( \pi x_1)  \sin (\pi x_2)  \\ 
0  \\ 0  \\ 0
\end{bmatrix}
\end{equation}
we have 
\begin{equation}
\begin{split}
 \|{\bf u}_A(x,0)\|^2_{L^2(\Omega)}&=({\bf u}_A(x,0),{\bf u}_A(x,0))\\
 &=\int_{(0,1)^3} ( \gamma^2 \sin^2 (\pi x_2)  \sin^2 (\pi x_3) + 4\gamma^2 \sin^2 (\pi x_1) \sin^2 ( \pi x_3)\\
 &+9 \gamma^2 \sin^2 ( \pi x_1)  \sin^2 (\pi x_2)  )dx_1 dx_2 dx_3=\frac{1}{4}\gamma^2  + \gamma^2  +   \frac{9}{4} \gamma^2  
 \end{split}
\end{equation}
since $\int_{(0,1)}\sin^2(\pi x_i)dx_i=\frac{1}{2}$. We want $\frac{14}{4} \gamma^2 =1$, 
so $\gamma=\sqrt{\frac{4}{14}}=\frac{2}{\sqrt{14}}.$

We summarize the numerical experiments in Figures 1-4. Figure 1 is the snapshot from the numerical simulation. Figure 2 shows that we have an implicit second-order in time method. Figures 3 and 4 illustrate how the $L^2$ and $H-curl$ error of the method changes when we increase the number of time steps, from $\tau=\frac{1}{10}$ down to $\tau=\frac{1}{1280}$. The $L^2$ error for $\tau=\frac{1}{10}$ is less than $0.08$, and for $\tau=\frac{1}{1280}$ it is less than 0.0002.
The $H-curl$ error for $\tau=0.1$ is less than $0.35$, and for $\tau=\frac{1}{1280}$ it is less than 0.015.

\begin{figure}
	\begin{center}
	\includegraphics[scale=0.2]{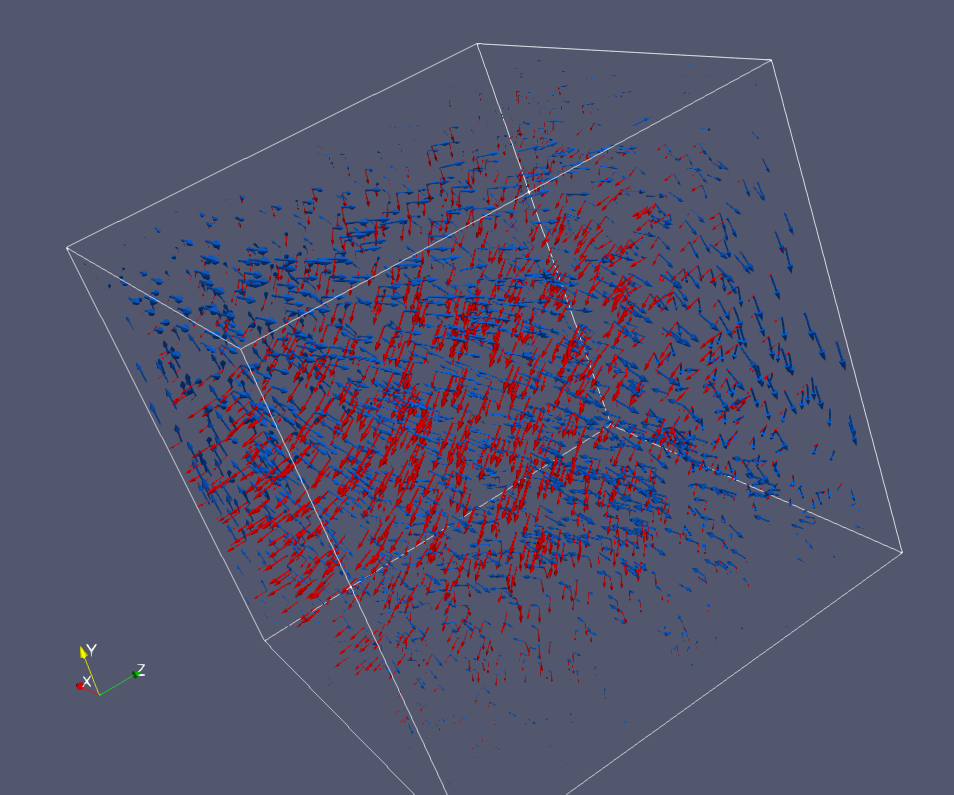}
	\end{center}
	\caption{Electric (red) and magnetic (blue) vector fields, resulting from the problem with manufactured solution.}
	\label{fig:Maxwell1}
\end{figure}

\begin{figure}
	\begin{center}
	\includegraphics[scale=0.8]{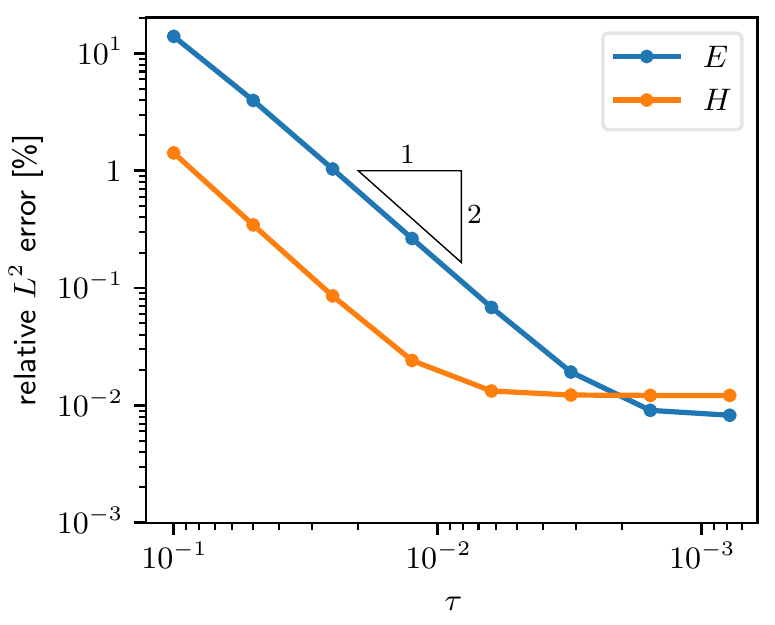}
	\includegraphics[scale=0.8]{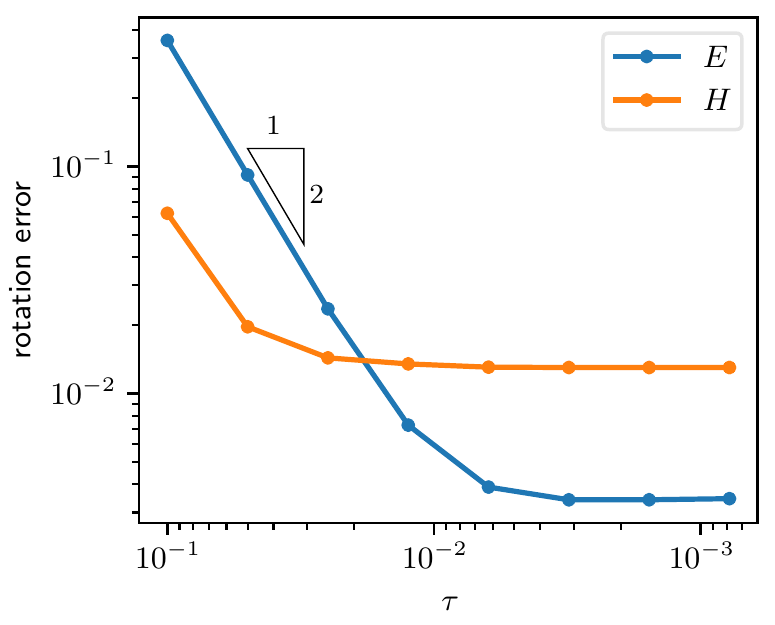}
	\end{center}
	\caption{Order of the time integration scheme as measured in L2 (left) and H-curl (right) norms for electric (blue) and magnetic (orange) vector fields resulting from the solution of the problem with manufactured solution over the computational mesh with 16x16x16 elements.}
	\label{fig:Maxwell1a}
\end{figure}

\begin{figure}
	\begin{center}
	\includegraphics[scale=0.6]{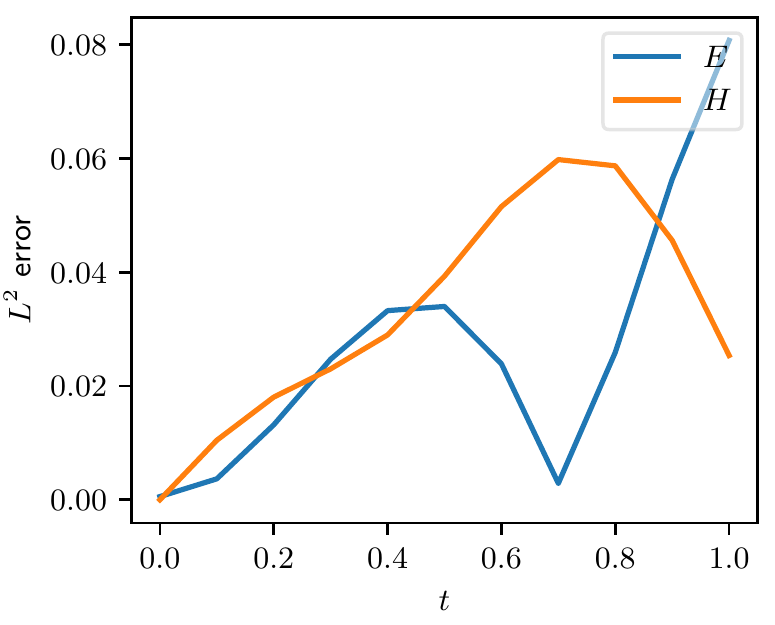}
	\includegraphics[scale=0.6]{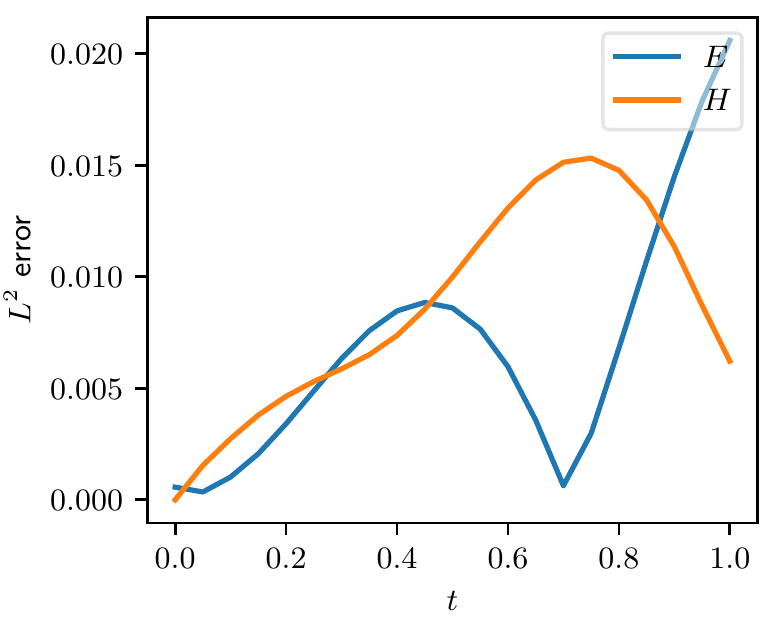}
	\includegraphics[scale=0.6]{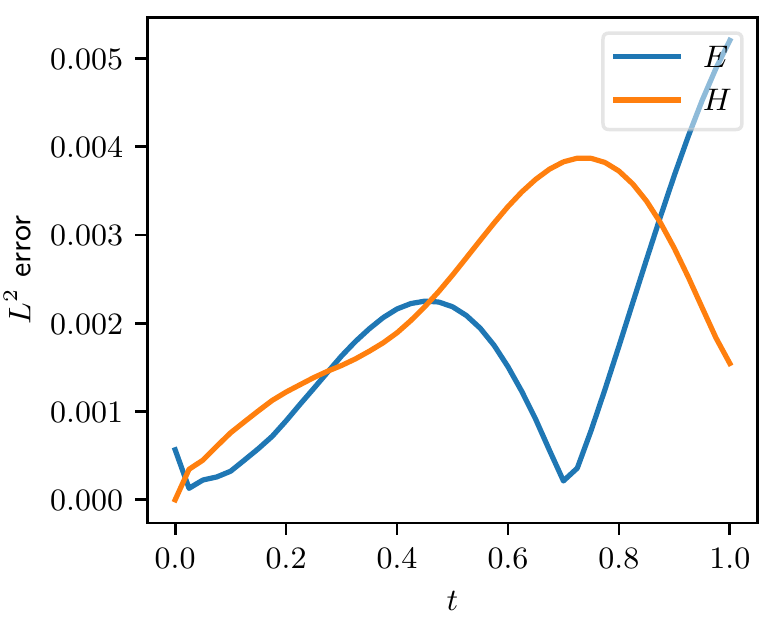}
	\includegraphics[scale=0.6]{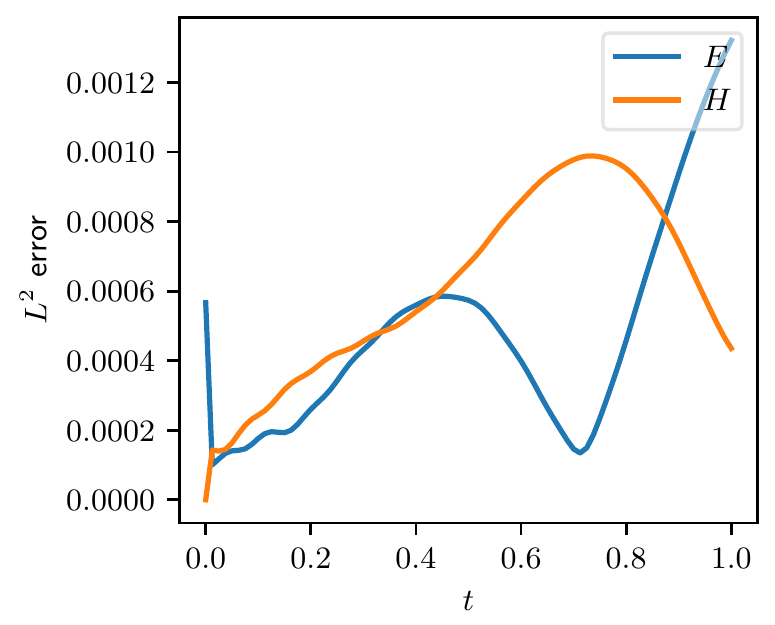}
	\includegraphics[scale=0.6]{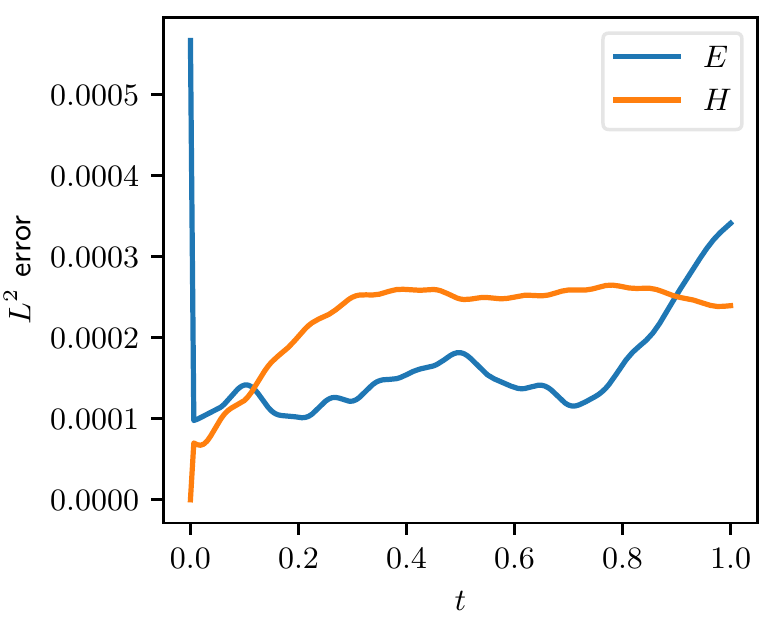}
	\includegraphics[scale=0.6]{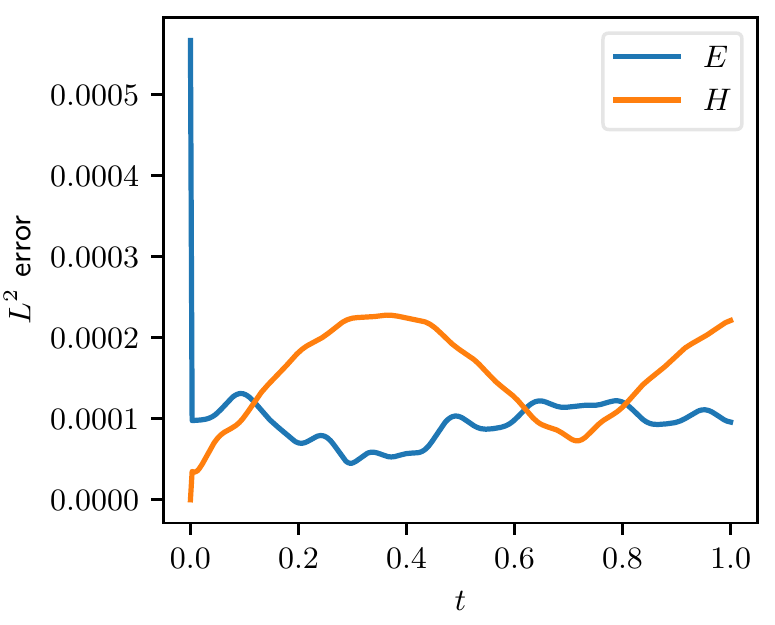}
	\includegraphics[scale=0.6]{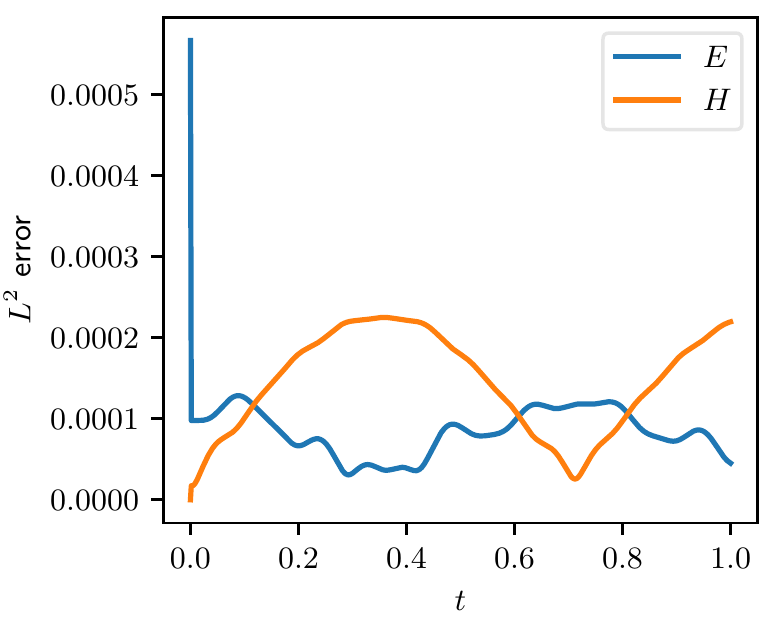}
	\includegraphics[scale=0.6]{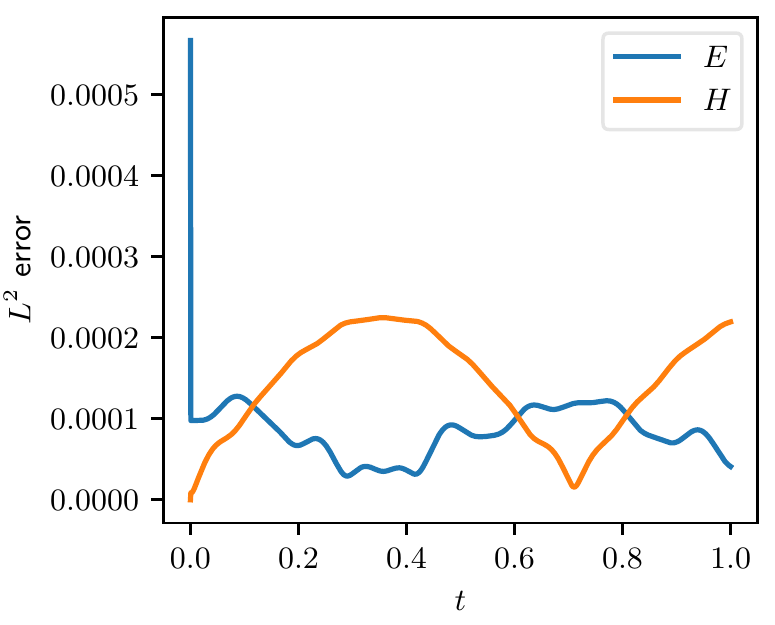}
	\end{center}
	\caption{L2 norm error of electric (blue) and magnetic (orange) vector fields resulting from the solution of the problem with manufactured solution over the computational mesh with 16x16x16 elements, for the time interval [0,1], with number of time steps within [0,1] interval varying from 10,20 (first row), 40,80 (second row), 160,320 (third row), 640 and 1280 (last row).}
	\label{fig:Maxwell1a}
\end{figure}
\begin{figure}
	\begin{center}
	\includegraphics[scale=0.6]{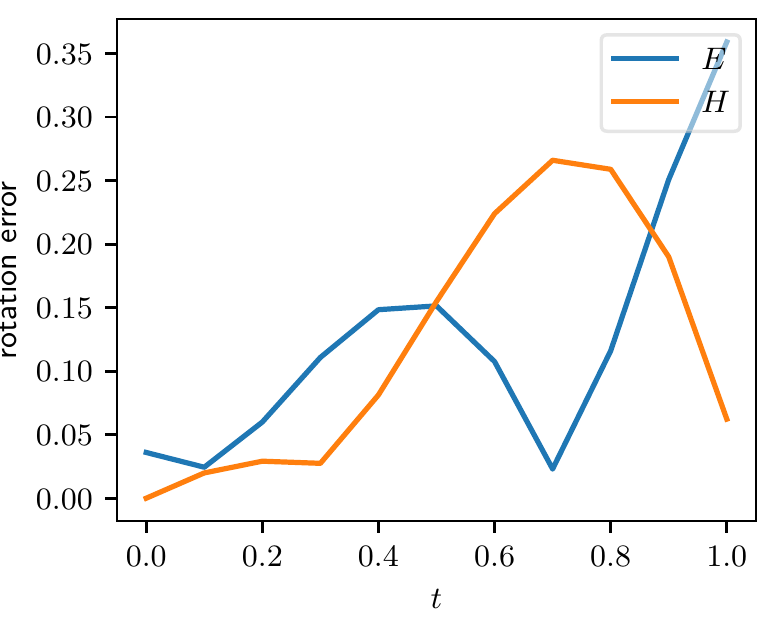}
	\includegraphics[scale=0.6]{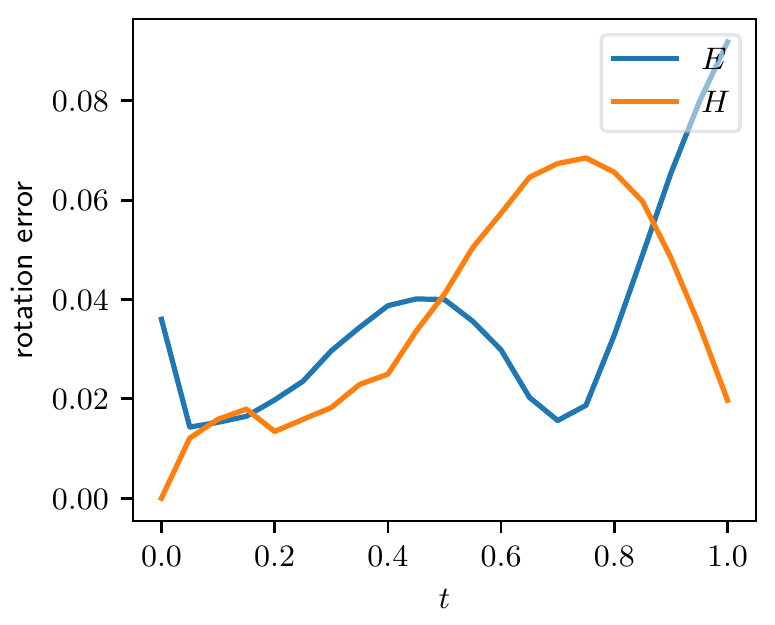}
	\includegraphics[scale=0.6]{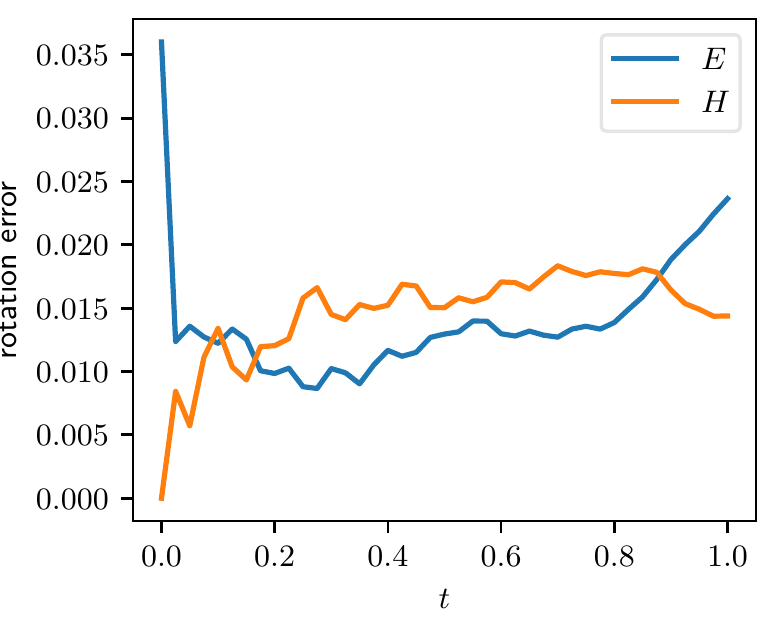}
	\includegraphics[scale=0.6]{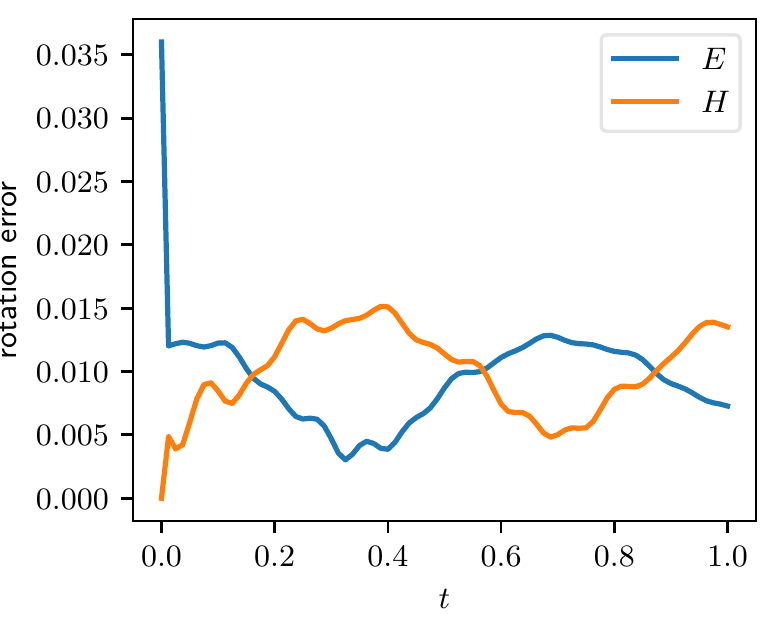}
	\includegraphics[scale=0.6]{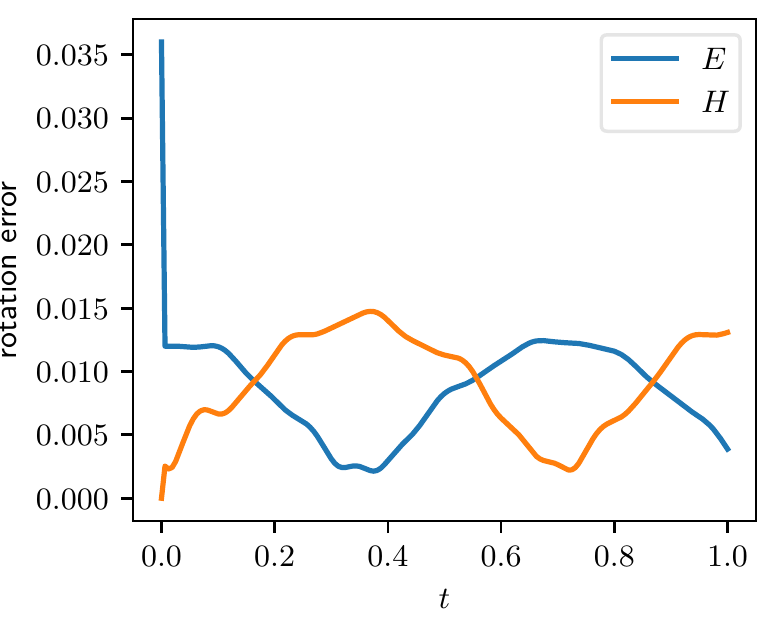}
	\includegraphics[scale=0.6]{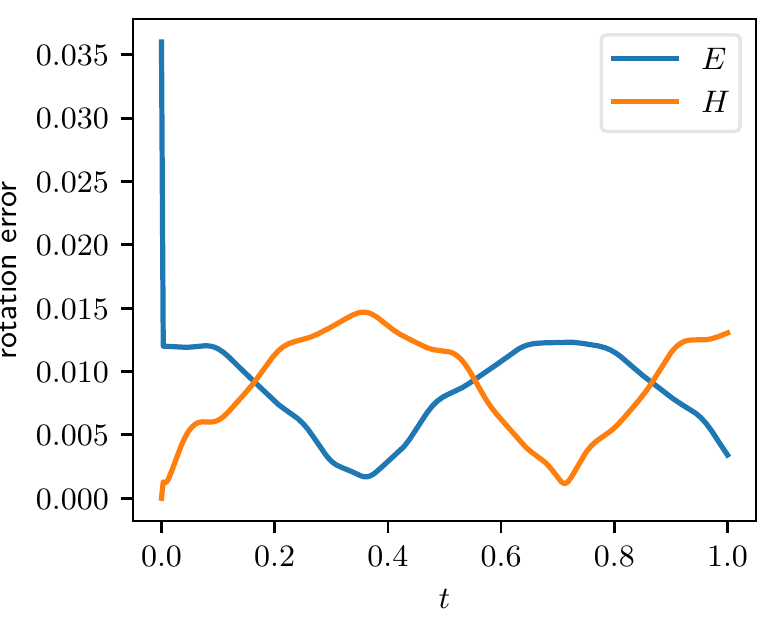}
	\includegraphics[scale=0.6]{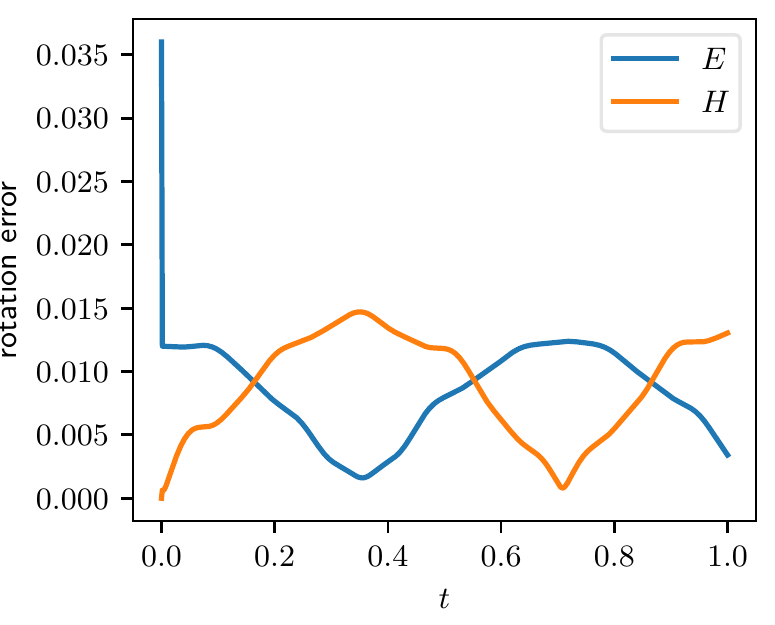}
	\includegraphics[scale=0.6]{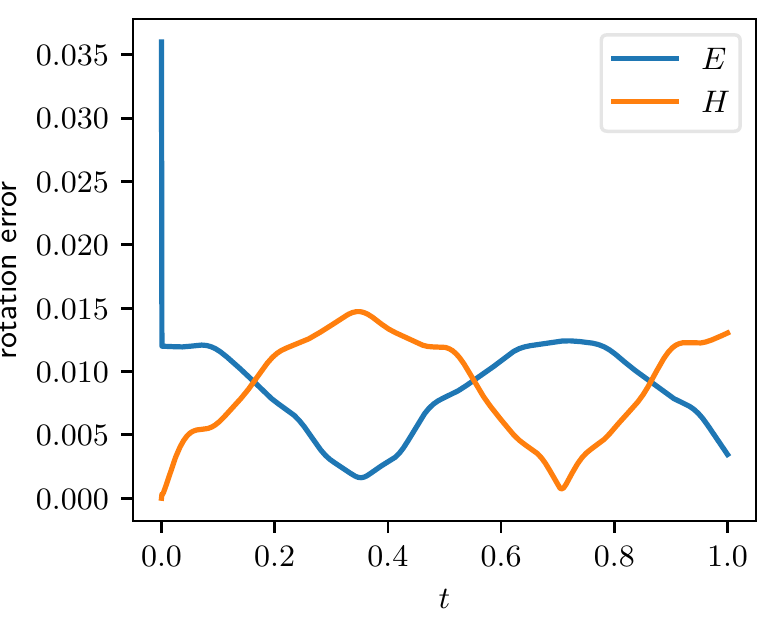}
	\end{center}
	\caption{H-curl norm error of electric (blue) and magnetic (orange) vector fields resulting from the solution of the problem with manufactured solution over the computational mesh with 16x16x16 elements, for the time interval [0,1], with number of time step varying from 10,20 (first row), 40,80 (second row), 160,320 (third row), 640 and 1280 (last row).}
	\label{fig:Maxwell1a}
\end{figure}

\section*{Incorporating non-regular material data into isogeometric alternating-direction solver}

We utilize the alternating directions solver that delivers linear computational cost factorization on tensor product grids. The solver decomposes the system of linear equations related to the three-dimensional mesh into three multi-diagonal sub-systems related to one-dimensional grids with multiple right-hand sides. The non-regular material data can be embedded into the solver by local modifications to the rows and columns in the three sub-systems. Namely, we can change the material data corresponding to different equations, and these modifications do not break the solver's linear computational cost. We verify this method by running the example of propagation of electromagnetic waves on the human head. Petar Minev has proposed this method initially for finite difference simulations \cite{Minev}. In the IGA context, the modification is not point-wise but rather test-function-wise since each equation in the global system is related to a single test function rather than a point in the stencil. 
Let us explain this idea in the example, using the first system of equations, solved in the even sub-steps, to update the electric field.  For other systems, the idea is identical. For simplicity in the notation we employ now $\{x,y,z\}$ instead of $\{x_1,x_2,x_3\}$. In the problem matrix, for the even sub-steps, for the electric field computations, we have  after multiplying the block matrices
\begin{equation}\label{dupa1ab}
\begin{split}
\begin{bmatrix}   M_{x} \otimes \left( M_{y} + \frac{\tau^2 }{4 \epsilon \mu}S_{y}\right)\otimes  M_{z} E_1^{n+\frac{1}{2}} \\ M_{x} \otimes M_{y} \otimes  \left(M_{z}+\frac{\tau^2 }{4 \epsilon \mu}S_{z}\right) E_2^{n+\frac{1}{2}} \\ \left(M_{x}+\frac{\tau^2 }{4 \epsilon \mu}S_{x}\right)\otimes M_{y}\otimes M_{z} E_3^{n+\frac{1}{2}} \end{bmatrix} \\ 
=\begin{bmatrix}  M_{x} \otimes M_{y}\otimes M_{z} E_1^n \\M_{x}\otimes  M_{y}\otimes M_{z} E_2^n \\ M_{x} \otimes M_{y}\otimes M_{z} E_3^n\end{bmatrix} +
 \begin{bmatrix}- \frac{\tau }{2 \epsilon} M_{x}\otimes M_{y} \otimes A_{z} H^n_2 \\
 \frac{\tau }{2 \epsilon} M_{x} \otimes M_{y}\otimes A_{z} H^n_1  \\  - \frac{\tau }{2 \epsilon} M_{x} \otimes A_{y} \otimes M_{z} H^n_1   
\end{bmatrix}\\ 
+ \begin{bmatrix} \frac{\tau }{2 \epsilon} M_{x} \otimes A_{y} \otimes M_{z} H^n_3  \\
 - \frac{\tau }{2 \epsilon} A_{x} \otimes M_{y} \otimes M_{z} H^n_3 \\  \frac{\tau }{2 \epsilon}A_{x} \otimes M_{y} \otimes M_{z} H^n_2 \\ 
\end{bmatrix} 
+   \begin{bmatrix} \frac{\tau^2 }{4\epsilon\mu} A_{x} \otimes B_{y} \otimes M_{z} E_2^n \\ 
\frac{\tau^2 }{4\epsilon\mu} M_{x} \otimes A_{y} \otimes B_{z} E^n_3 \\  \frac{\tau^2 }{4\epsilon\mu} B_{x} \otimes M_{y} \otimes A_{z} E_1^n  
\end{bmatrix}
\end{split}
\end{equation}
Rewriting the equations in matrix form with the B-spline functions for trial and testing, we have

\begin{equation}
\begin{bmatrix}
{\cal M}_1^1 {E_1^{n+\frac{1}{2}}} \\
{\cal M}_2^1 {E_2^{n+\frac{1}{2}}} \\
{\cal M}_3^1 {E_3^{n+\frac{1}{2}}}
\end{bmatrix}
 = 
\begin{bmatrix}
{\cal M}E_1^n \\
{\cal M}E_2^n \\
{\cal M} E_3^n 
\end{bmatrix}
+
\begin{bmatrix}
{\cal F}_1^1 H_2^n \\
{\cal F}_2^1 H_1^n \\
{\cal F}_3^1 H_1^n 
\end{bmatrix}
+
\begin{bmatrix}
{\cal F}_1^2 H_3^n \\
{\cal F}_2^2 H_3^n \\
{\cal F}_3^2 H_2^n
\end{bmatrix}
+
\begin{bmatrix}
{\cal F}_1^3 E_2^n\\
{\cal F}_2^3 E_3^n\\
{\cal F}_3^3 E_1^n
\end{bmatrix}
=\begin{bmatrix}
{\cal RHS}_1 \\
{\cal RHS}_2 \\
{\cal RHS}_3 \\
\end{bmatrix}
\end{equation}
where the entries of each matrix are
\begin{equation}
\begin{split}
{{\cal M}_1^1}_{ijk,lmo}&=\int_{\Omega_x} B_{i,p}(x) B_{l,p}(x) dx\\
& \int_{\Omega_y} \left( B_{j,p}(y)B_{m,p}(y) + \frac{\tau^2 }{4 \epsilon \mu}\frac{\partial B_{j,p}(y)}{\partial y}\frac{\partial B_{m,p}(y)}{\partial y}\right) dy\int_{\Omega_z} B_{k,p}(z)B_{o,p}(z)dz 
\\ 
{{\cal M}_2^1}_{ijk,lmo}&=\int_{\Omega_x} B_{i,p}(x) B_{l,p}(x) dx\int_{\Omega_y}  B_{j,p}(y) B_{m,p}(y)dy \\
&\int_{\Omega_z} \left( B_{k,p}(z)B_{n,p}(z) + \frac{\tau^2 }{4 \epsilon \mu}\frac{\partial B_{k,p}(z)}{\partial z}\frac{\partial B_{o,p}(z)}{\partial z}\right)dz 
\\
{{\cal M}_3^1}_{ijk,lmo}&=\int_{\Omega_x} \left( B_{i,p}(x)B_{l,p}(x) + \frac{\tau^2 }{4 \epsilon \mu}\frac{\partial B_{i,p}(x)}{\partial x}\frac{\partial B_{l,p}(x)}{\partial x}\right)dx \\
&\int_{\Omega_y}  B_{j,p}(y) B_{m,p}(y)dy\int_{\Omega_z} B_{k,p}(z)B_{o,p}(z)dz \\
\end{split}
\end{equation}
\begin{equation}
{{\cal M}}_{ijk,lmo}=\int_{\Omega_x} B_{i,p}(x) B_{l,p}(x) dx\int_{\Omega_y}B_{j,p}(y)B_{m,p}(y)dy\int_{\Omega_z} B_{k,p}(z)B_{o,p}(z)dz 
\end{equation}
\begin{equation}\label{dupa1abcd}
\begin{split}
{{\cal F}_1^1}_{ijk,lmo}&=
 -\frac{\tau }{2 \epsilon}
  \int_{\Omega} B_{i,p}(x) B_{j,p}(y) \frac{\partial B_{k,p}(z)}{\partial z} B_{l,p}(x)  B_{m,p}(y) B_{o,p}(z) dxdydz\\
{{\cal F}_2^1}_{ijk,lmo}&=
\frac{\tau }{2 \epsilon}\int_{\Omega} B_{i,p}(x) B_{j,p}(y) \frac{\partial B_{k,p}(z)}{\partial z} B_{l,p}(x)  B_{m,p}(y) B_{o,p}(z) dxdydz\\
{{\cal F}_3^1}_{ijk,lmo}&=
 -\frac{\tau }{2 \epsilon} \int_{\Omega} B_{i,p}(x) \frac{\partial B_{j,p}(y)}{\partial y}B_{k,p}(x) B_{l,p}(x)  B_{m,p}(y) B_{o,p}(z) dxdydz\\
{{\cal F}_1^2}_{ijk,lmo}&=
 \frac{\tau }{2 \epsilon} \int_{\Omega} B_{i,p}(x) \frac{\partial B_{j,p}(y)}{\partial y}B_{k,p}(x) B_{l,p}(x)  B_{m,p}(y) B_{o,p}(z) dxdydz\\
{{\cal F}_2^2}_{ijk,lmo}&=
 -\frac{\tau }{2 \epsilon}\int_{\Omega} \frac{\partial B_{i,p}(x)}{\partial x} B_{j,p}(y)B_{k,p}(x) B_{l,p}(x)  B_{m,p}(y) B_{o,p}(z) dxdydz\\  
{{\cal F}_3^2}_{ijk,lmo}&=
- \frac{\tau }{2 \epsilon} \int_{\Omega} \frac{\partial B_{i,p}(x)}{\partial x} B_{j,p}(y)B_{k,p}(x) B_{l,p}(x)  B_{m,p}(y) B_{o,p}(z) dxdydz\\
{{\cal F}_1^3}_{ijk,lmo}&=
\frac{\tau^2 }{4\epsilon\mu} \int_{\Omega} \frac{\partial B_{i,p}(x)}{\partial x}B_{k,p}(x)B_{j,p}(y)B_{l,p}(x)   \frac{\partial B_{m,p}(y)}{\partial y}  B_{o,p}(z) dxdydz\\
{{\cal F}_2^3}_{ijk,lmo}&=
\frac{\tau^2 }{4\epsilon\mu} \int_{\Omega}  B_{i,p}(x)\frac{\partial B_{j,p}(y)}{\partial y} B_{k,p}(z) B_{l,p}(x)  B_{m,p}(y) \frac{\partial B_{o,p}(z)}{\partial z} dxdydz\\
{{\cal F}_3^3}_{ijk,lmo}&=
\frac{\tau^2 }{4\epsilon\mu}\int_{\Omega} B_{i,p}(x)B_{j,p}(y)\frac{\partial B_{k,p}(z)}{\partial z}\frac{\partial B_{l,p}(x)}{\partial x} B_{m,p}(y) B_{o,p}(z) dxdydz 
\end{split}
\end{equation}
where $i=1,...,N_x$, $j=1,...,N_y$, $k=1,...,N_z$ span over the trial space dimensions, and $l=1,...,\tilde{N}_x$, $m=1,...,\tilde{N}_y$, $n=1,...,\tilde{N}_z$ span over the test space dimensions.
The matrices on the right-hand side are multiplied by the solution vectors from previous time step, so as the result on the right-hand side we have a vectors ${{\cal RHS}_1}_{lmo}$, ${{\cal RHS}_2}_{lmo}$, and ${{\cal RHS}_3}_{lmo}$, where again $l=1,...,\tilde{N}_x$, $m=1,...,\tilde{N}_y$, $o=1,...\tilde{N}_z$.

The alternating directions solver decomposes this system into the following three one-dimensional systems with multiple right-hand-sides
\begin{eqnarray}
\begin{bmatrix}   
{\cal A}_1 {F_1^{n+\frac{1}{2}}} \\
{\cal A}_2 {F_2^{n+\frac{1}{2}}} \\
{\cal A}_3 {F_3^{n+\frac{1}{2}}} \\
\end{bmatrix}   
=
\begin{bmatrix}
{\cal RHS}_1 \\
{\cal RHS}_2 \\
{\cal RHS}_3 \\
\end{bmatrix}
\end{eqnarray}
where
\begin{eqnarray}
{{\cal A}_1}_{i,l} = 
\int_{\Omega_x} \inred{B_{i,p}(x)}\inblue{ B_{l,p}(x)} dx
\nonumber \\ 
{{\cal A}_2}_{i,l} = 
  \int_{\Omega_x} \inred{B_{i,p}(x)} \inblue{ B_{l,p}(x)} dx  
 \nonumber \\ 
{{\cal A}_3}_{i,l} = 
\int_{\Omega_x} \left( \inred{B_{i,p}(x)}\inblue{B_{l,p}(x)} + \frac{\tau^2 }{4 \epsilon \mu}\inred{\frac{\partial B_{i,p}(x)}{\partial x}}\inblue{\frac{\partial B_{l,p}(x)}{\partial x}}\right)dx 
\end{eqnarray}
and the right-hand side vectors ${{\cal RHS}_1}_{i,jk}$, ${{\cal RHS}_2}_{i,jk}$, ${{\cal RHS}_3}_{i,jk}$ have been reordered into matrices with $N_x$ rows and $N_y N_z$ columns, 
by ordering blocks of $N_x$ consecutive rows, one after another.

After solving the first one-dimensional system with multiple right-hand sides we solve the second system
\begin{eqnarray}
\begin{bmatrix}   
{\cal B}_1G_1^{n+\frac{1}{2}}  \\
{\cal B}_2G_2^{n+\frac{1}{2}} \\
{\cal B}_3G_3^{n+\frac{1}{2}} \\
\end{bmatrix}   
=
\begin{bmatrix}F_1^{n+\frac{1}{2}} \\ F_2^{n+\frac{1}{2}} \\ F_3^{n+\frac{1}{2}} \end{bmatrix} 
\end{eqnarray}
where
\begin{eqnarray}
{{\cal B}_1}_{j,m}= \int_{\Omega_y} \left( \inred{B_{j,p}(y)}\inblue{B_{m,p}(y)} + \frac{\tau^2 }{4 \epsilon \mu}\inred{\frac{\partial B_{j,p}(y)}{\partial y}}\inblue{\frac{\partial B_{m,p}(y)}{\partial y}}\right) dy \nonumber \\
{{\cal B}_2}_{j,m}=  \int_{\Omega_y}  \inred{B_{j,p}(y)}\inblue{ B_{m,p}(y)}dy \nonumber \\ 
{{\cal B}_3}_{j,m}= \int_{\Omega_y}  \inred{B_{j,p}(y)}\inblue{ B_{m,p}(y)}dy
\end{eqnarray}
Finally, we solve the third system with multiple right-hand sides
\begin{eqnarray}
\begin{bmatrix}   
{\cal C}_1 E_1^{n+\frac{1}{2}} \\
{\cal C}_2 E_2^{n+\frac{1}{2}} \\
{\cal C}_3 E_3^{n+\frac{1}{2}} \\
\end{bmatrix}   
=
\begin{bmatrix}G_1^{n+\frac{1}{2}} \\ G_2^{n+\frac{1}{2}} \\ G_3^{n+\frac{1}{2}} \end{bmatrix} 
\end{eqnarray}
where
\begin{eqnarray}
{{\cal C}_1}_{k,o}= \int_{\Omega_z} \inred{B_{k,p}(z)}\inblue{B_{o,p}(z)}dz \nonumber \\ 
{{\cal C}_2}_{k,o}=\int_{\Omega_z} \left( \inred{B_{k,p}(z)}\inblue{B_{o,p}(z)} + \frac{\tau^2 }{4 \epsilon \mu}\inred{\frac{\partial B_{k,p}(z)}{\partial z}}\inblue{\frac{\partial B_{o,p}(z)}{\partial z}}\right)dz \nonumber \\ 
{{\cal C}_3}_{k,o}=\int_{\Omega_z} \inred{B_{k,p}(z)}\inblue{B_{o,p}(z)}dz 
\end{eqnarray}

We need to modify the material data of the Maxwell equations related to tissue, skull, and air. We assign different material data to different B-splines used for testing our equation. Since each test B-spline results in a single equation in the global system of equations, we identify this equation in the three systems with multiple right-hand sides. Having the equations identified, we modify the material data in the three systems of equations as processed by the alternating directions solver.

For example, if we want to modify material data $\tau = \hat{\tau}$, $\epsilon=\hat{\epsilon}$, $\mu=\hat{\mu}$ for test B-spline "$rst$", namely \inblue{$B_{r,p}(x)B_{s,p}(y)B_{t,p}(z)$} we perform the following changes.
In the first system, we extract the three equations (three rows) for the three components of the electric field for row $i=r$, 
and the suitable columns from the right-hand side $l=r,m=s,o=t$, where we modify the material data

\begin{eqnarray}
\sum_{l=1,...,N_x} \int_{\Omega_x} \inred{B_{r,p}(x)}\inblue{ B_{l,p}(x)} dx 
{F_1^{n+\frac{1}{2}}}_{lst} = \hat{\cal RHS}_{1rst} 
\label{dupa1abcdefg}
\end{eqnarray}

\begin{eqnarray}
\sum_{l=1,...,N_x}\int_{\Omega_x} \inred{B_{r,p}(x)} \inblue{ B_{l,p}(x)} dx   
  {F_2^{n+\frac{1}{2}}}_{lst} = \hat{\cal RHS}_{2rst} \end{eqnarray}
\begin{eqnarray}
\sum_{l=1,...,N_x}\int_{\Omega_x} \left( \inred{B_{r,p}(x)}\inblue{B_{l,p}(x)} + \frac{\hat{\tau}^2 }{4 \hat{\epsilon} \hat{\mu}}\inred{\frac{\partial B_{r,p}(x)}{\partial x}}\inblue{\frac{\partial B_{l,p}(x)}{\partial x}}\right)dx 
{F_3^{n+\frac{1}{2}}}_{lst} = \hat{\cal RHS}_{3rst} 
\end{eqnarray}
The $\hat{\cal RHS}_{1rst}$, $\hat{\cal RHS}_{2rst}$,
$\hat{\cal RHS}_{3rst}$ represent the right-hand sides with material data parameters $\tau = \hat{\tau}$, $\epsilon=\hat{\epsilon}$, $\mu=\hat{\mu}$.
The other rows and columns in the first system remain unchanged.

Similarly, in the second system, we extract the equation for row $j=s$ and columns
$l=r,m=s,n=t$
\begin{eqnarray}
\sum_{m=1,...,N_y}\int_{\Omega_y} \left( \inred{B_{s,p}(y)}\inblue{B_{m,p}(y)} + \frac{\hat{\tau}^2 }{4 \hat{\epsilon} \hat{\mu}}\inred{\frac{\partial B_{s,p}(y)}{\partial y}}\inblue{\frac{\partial B_{m,p}(y)}{\partial y}}\right) dy 
{G_1^{n+\frac{1}{2}}}_{rmt} =
{F_1^{n+\frac{1}{2}}}_{rst} 
\end{eqnarray}
\begin{eqnarray}
\sum_{m=1,...,N_y}\int_{\Omega_y}  \inred{B_{s,p}(y)}\inblue{ B_{m,p}(y)dy} 
{G_2^{n+\frac{1}{2}}}_{rmt}   =
{F_2^{n+\frac{1}{2}}}_{rst} 
\end{eqnarray}
\begin{eqnarray}
\sum_{m=1,...,N_y} \int_{\Omega_y}  \inred{B_{s,p}(y)}\inblue{ B_{m,p}(y)dy}
{G_3^{n+\frac{1}{2}}}_{rmt}  =
{F_3^{n+\frac{1}{2}}}_{rst} 
\end{eqnarray}
and we modify the material data.
The other rows and columns in the second system remain unchanged.

Finally, in the third system, we extract the equation for row $k=t$ and columns
$l=r,m=s,n=t$
\begin{eqnarray}
\sum_{o=1,...,N_z}\int_{\Omega_z} \inred{B_{t,p}(z)}\inblue{B_{o,p}(z)}dz 
{E_1^{n+\frac{1}{2}}}_{rso}={G_1^{n+\frac{1}{2}}}_{rst} 
\end{eqnarray}

\begin{eqnarray}
\sum_{o=1,...,N_z}
\int_{\Omega_z} \left( \inred{B_{t,p}(z)}\inblue{B_{o,p}(z)} + \frac{\hat{\tau}^2 }{4 \hat{\epsilon} \hat{mu}}\inred{\frac{\partial B_{t,p}(z)}{\partial z}}\inblue{\frac{\partial B_{o,p}(z)}{\partial z}}\right)dz 
  {E_3^{n+\frac{1}{2}}}_{rso} =
{G_2^{n+\frac{1}{2}}}_{rst}
\end{eqnarray}

\begin{eqnarray}
\sum_{o=1,...,N_z}\int_{\Omega_z} \inred{B_{o,p}(z)}\inblue{B_{o,p}(z)}dz
 {E_3^{n+\frac{1}{2}}}_{rso} =
{G_3^{n+\frac{1}{2}}}_{rst} %\nonumber \\
%o=1,...,N_z
\end{eqnarray}
and we modify the material data.
The other rows and columns in the third system remain unchanged.

Similar modifications have to be performed in other sub-steps.

\begin{figure}
	\begin{center}
	\includegraphics[scale=0.4]{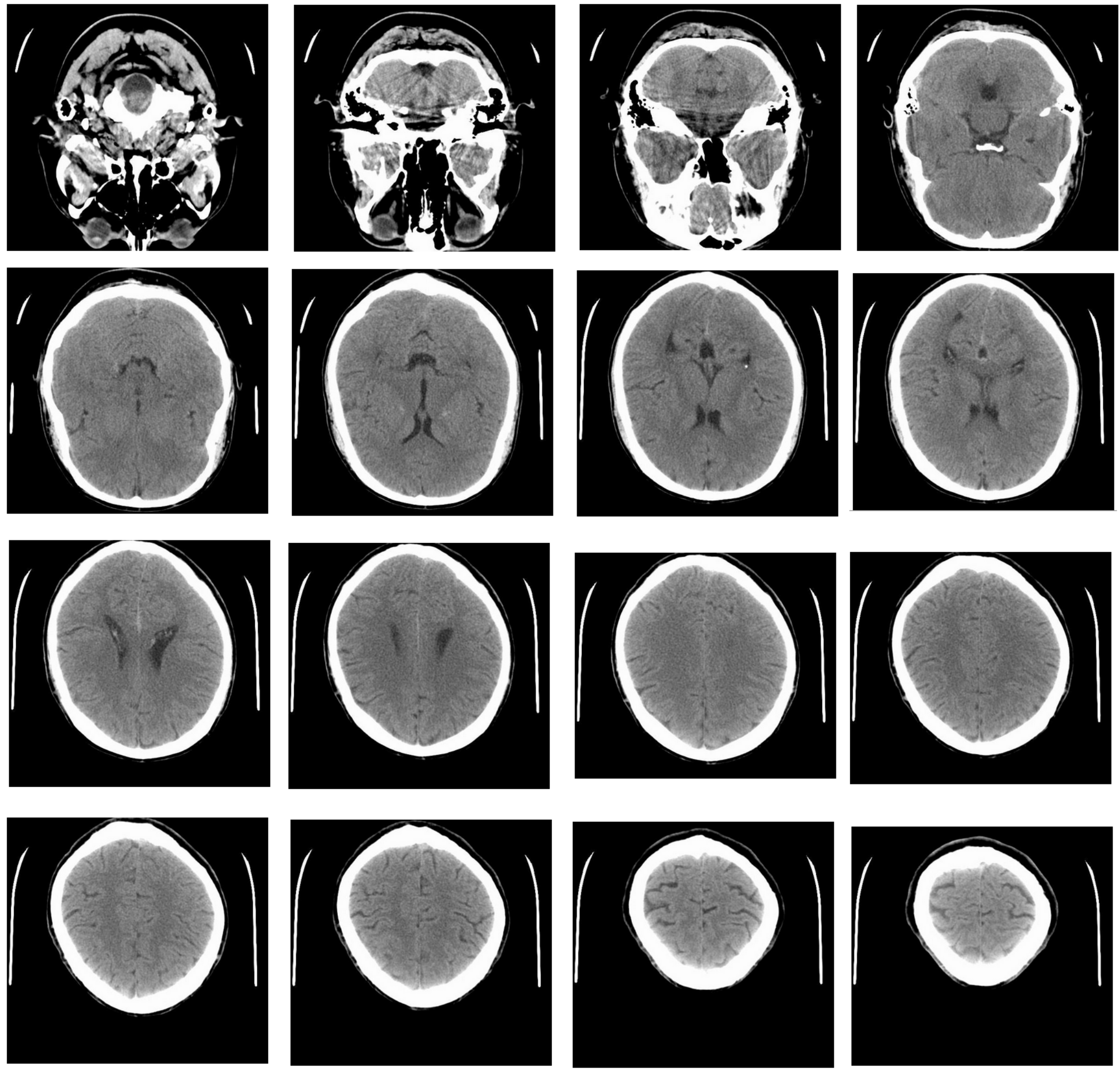}
	\end{center}
	\caption{Exemplary cross-sections of the MRI scans of the head of the first author.}
	\label{fig:head}
\end{figure}

\begin{figure}
	\begin{center}
	\includegraphics[scale=0.24]{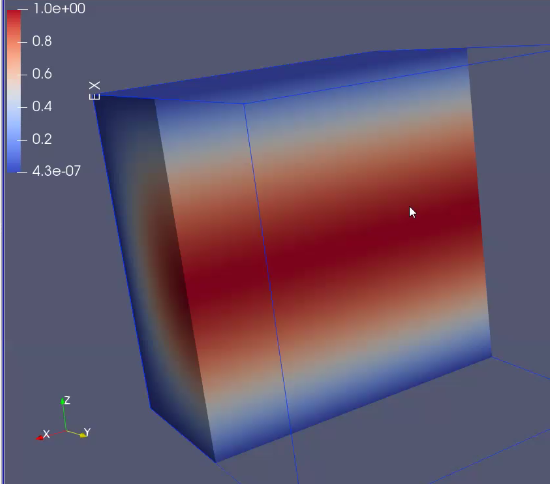}
	\includegraphics[scale=0.24]{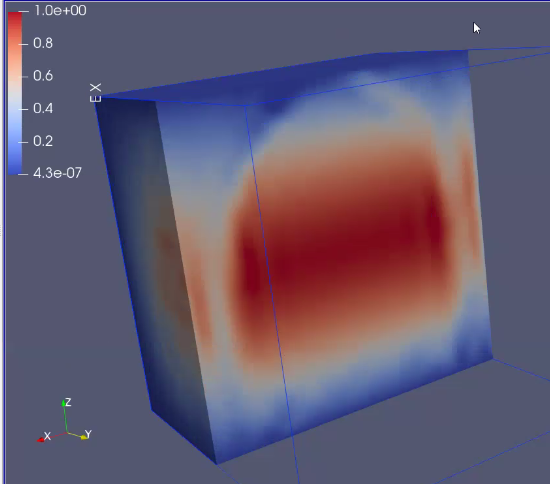}
	\includegraphics[scale=0.24]{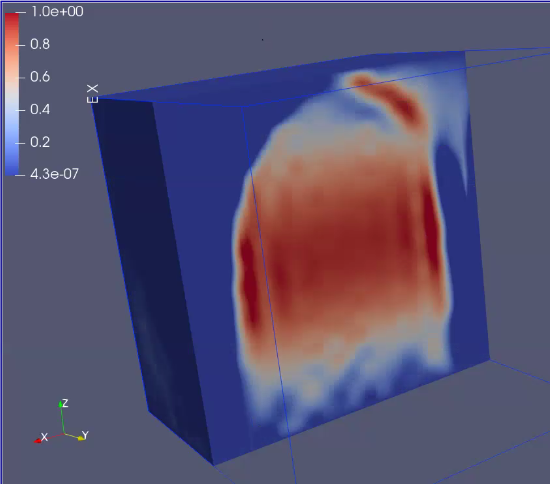}
	\includegraphics[scale=0.24]{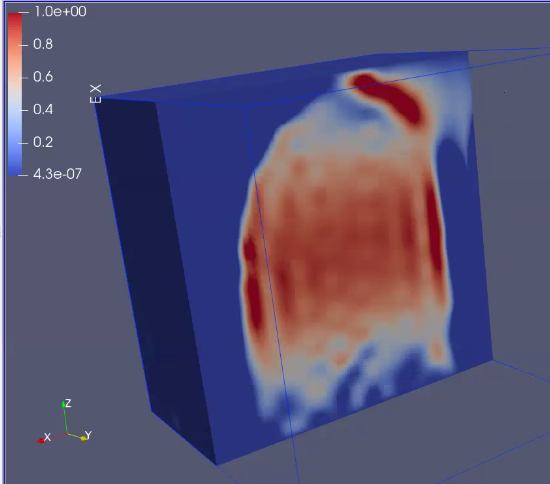}
	\includegraphics[scale=0.24]{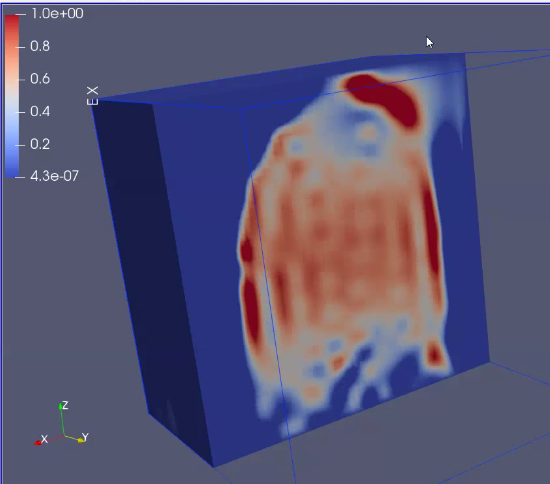}\\
	\includegraphics[scale=0.24]{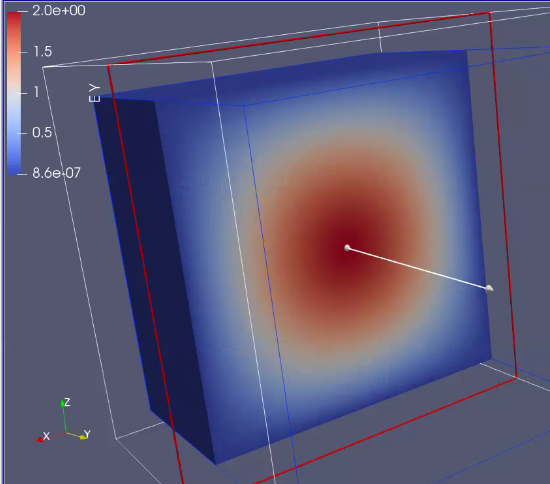}
	\includegraphics[scale=0.24]{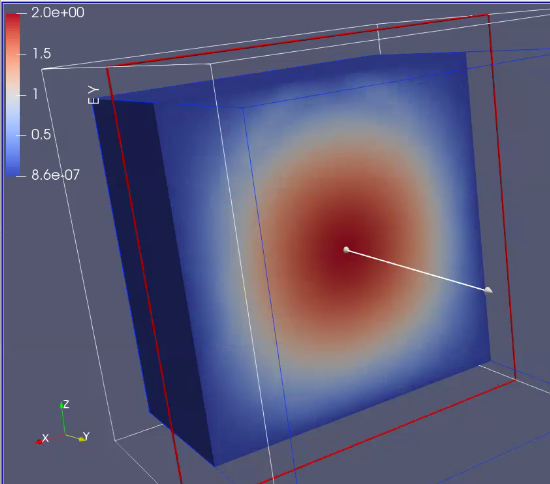}
	\includegraphics[scale=0.24]{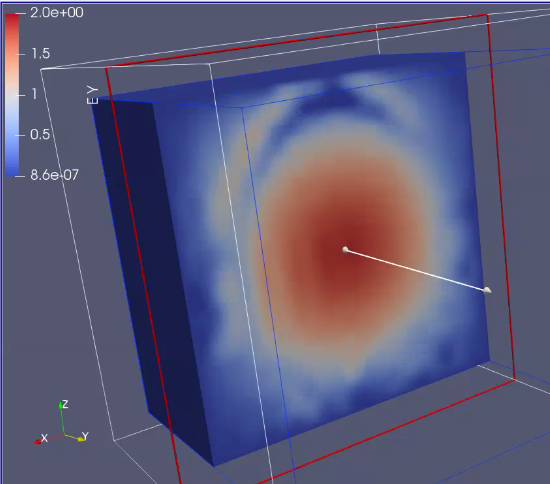}
	\includegraphics[scale=0.24]{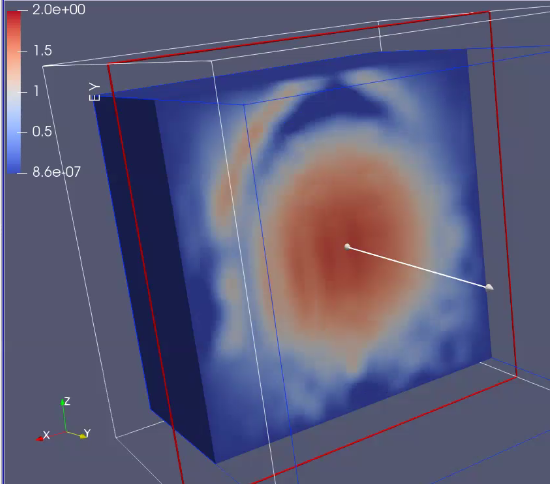}
	\includegraphics[scale=0.24]{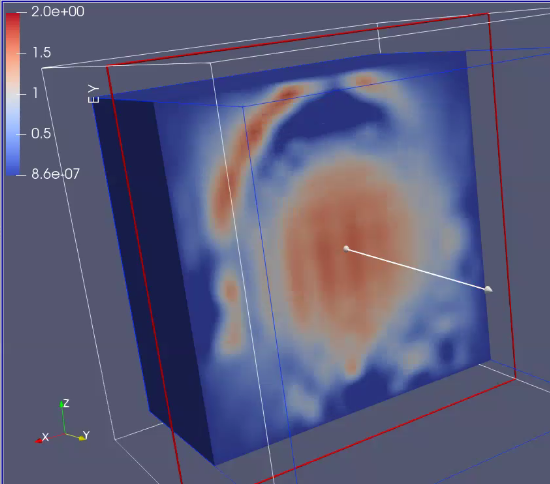}\\
	\includegraphics[scale=0.24]{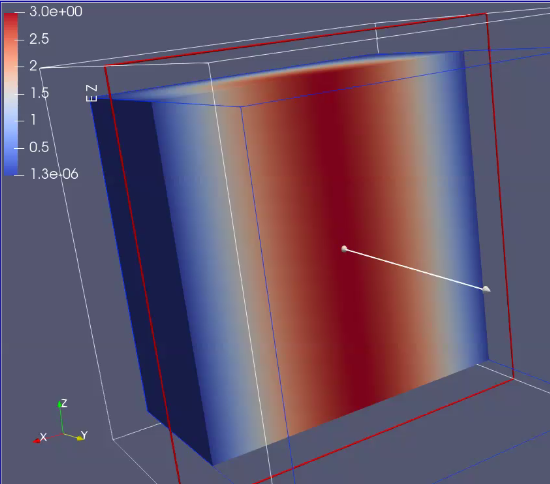}
	\includegraphics[scale=0.24]{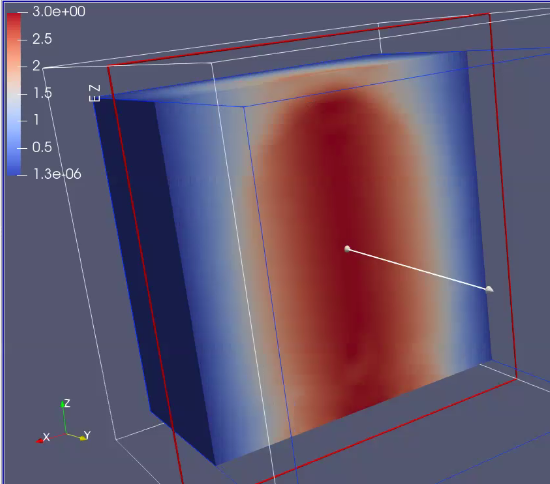}
	\includegraphics[scale=0.24]{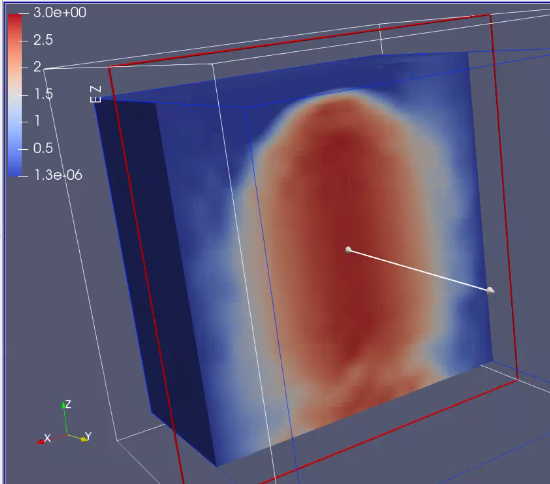}
	\includegraphics[scale=0.24]{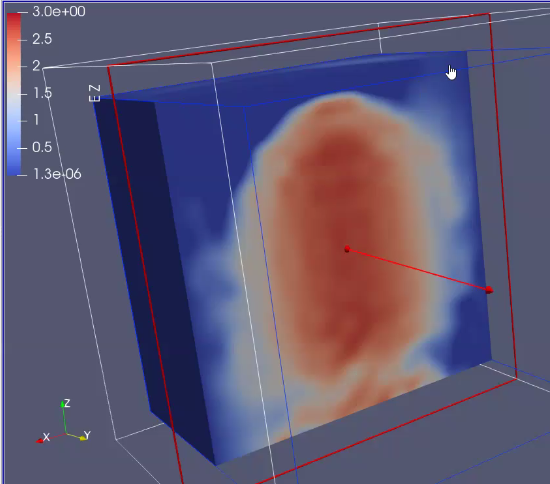}
	\includegraphics[scale=0.24]{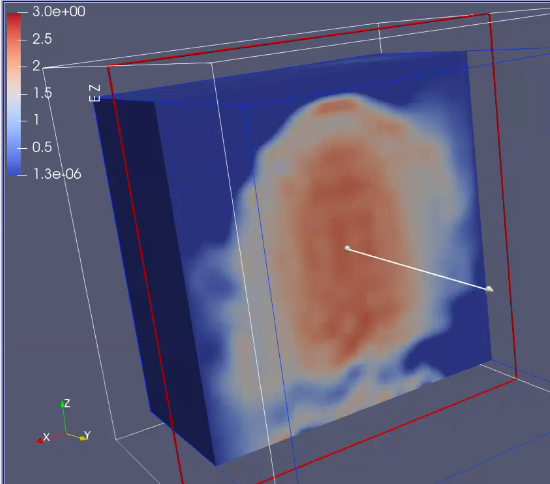}\\
	\end{center}
	\caption{Electromagnetic waves propagation on the human head. First row: $x$ component of the electric vector field. Second row: $y$ component of the electric vector field. Third row: $z$ component of the electric vector field. Columns - cross section along OYZ plane, time moments: 0, 0.25, 0.5, 0.75, 1.0s. Mesh size 32x32x32, quadratic $C^1$ B-splines.}
	\label{fig:head2}
\end{figure}

\begin{figure}
	\begin{center}
	\includegraphics[scale=0.2]{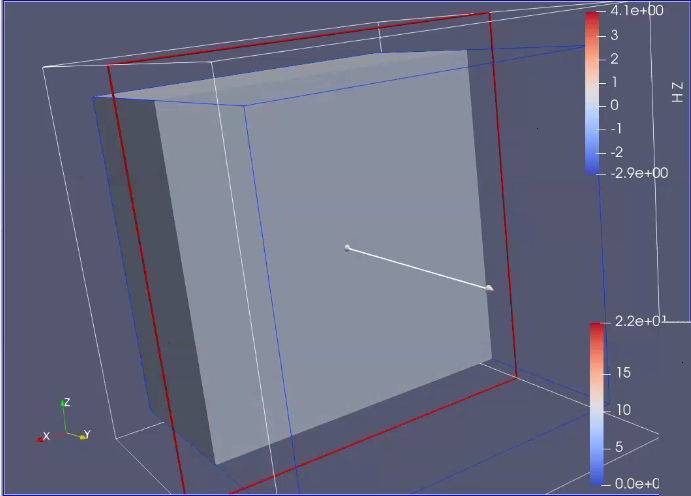}
	\includegraphics[scale=0.2]{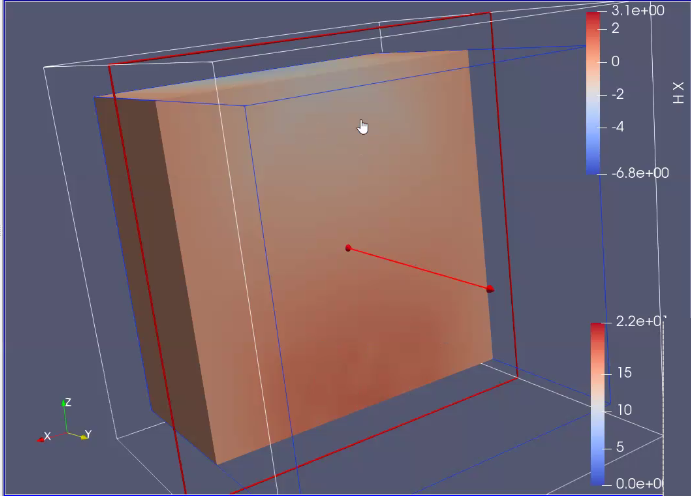}
	\includegraphics[scale=0.2]{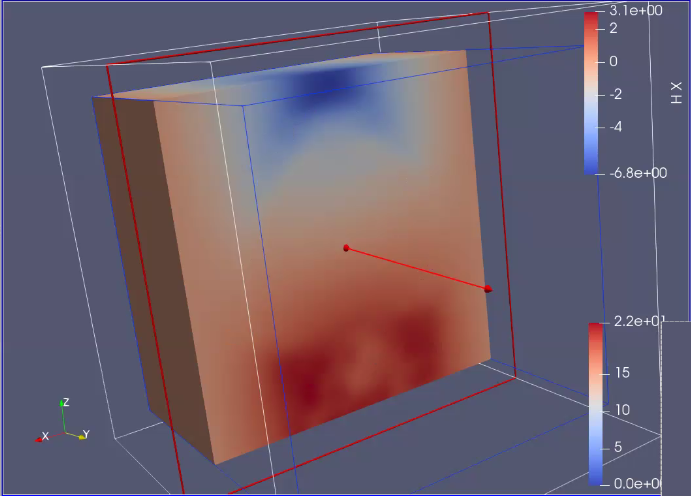}
	\includegraphics[scale=0.2]{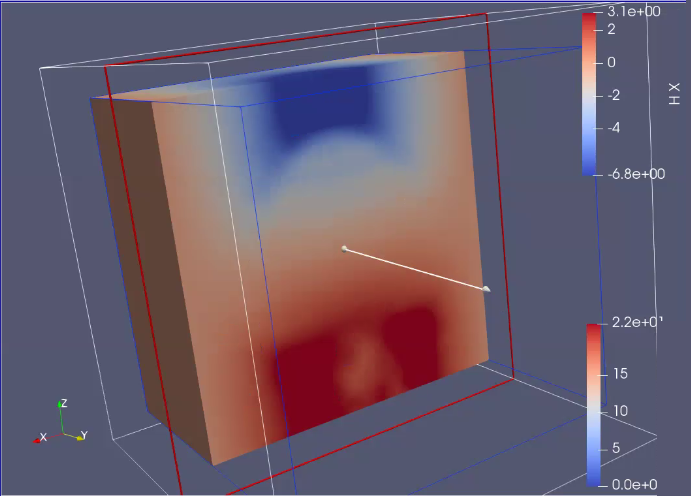}
	\includegraphics[scale=0.2]{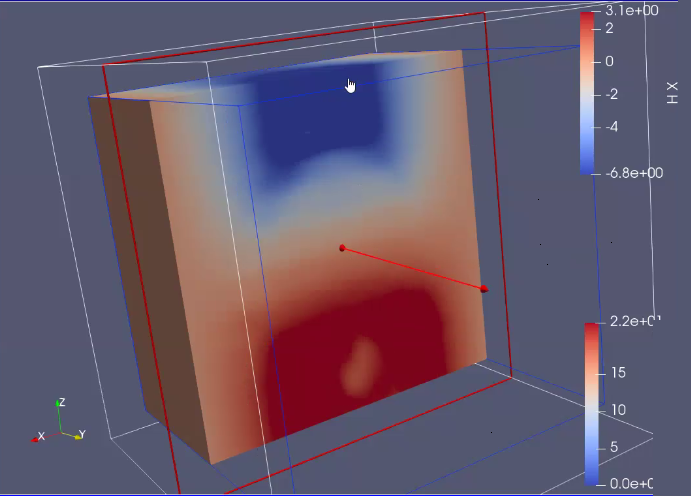}\\
	\includegraphics[scale=0.2]{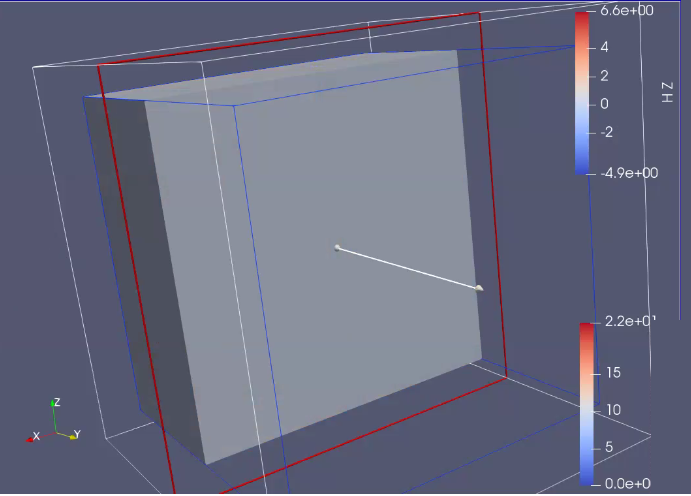}
	\includegraphics[scale=0.2]{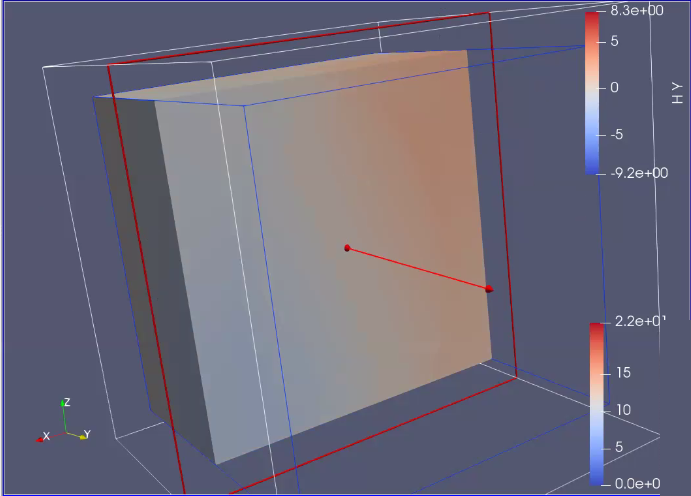}
	\includegraphics[scale=0.2]{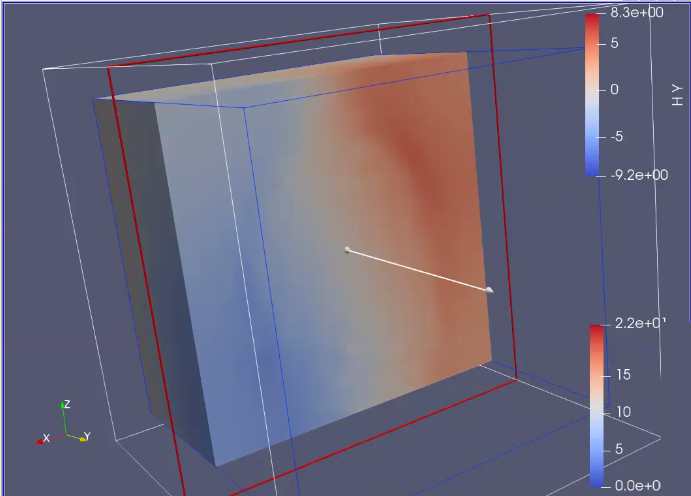}
	\includegraphics[scale=0.2]{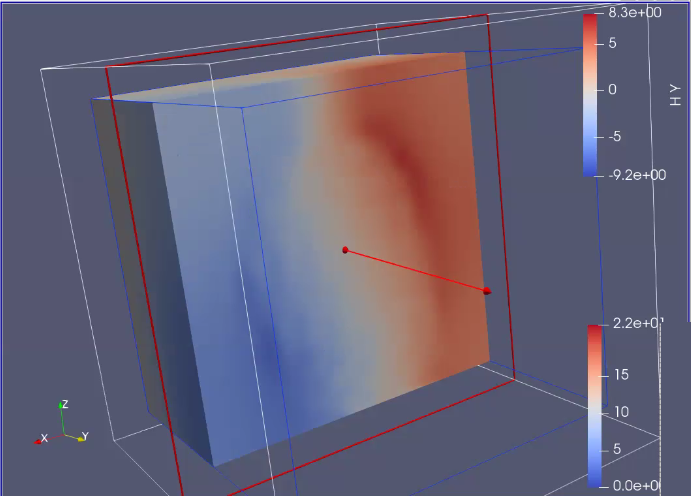}
	\includegraphics[scale=0.2]{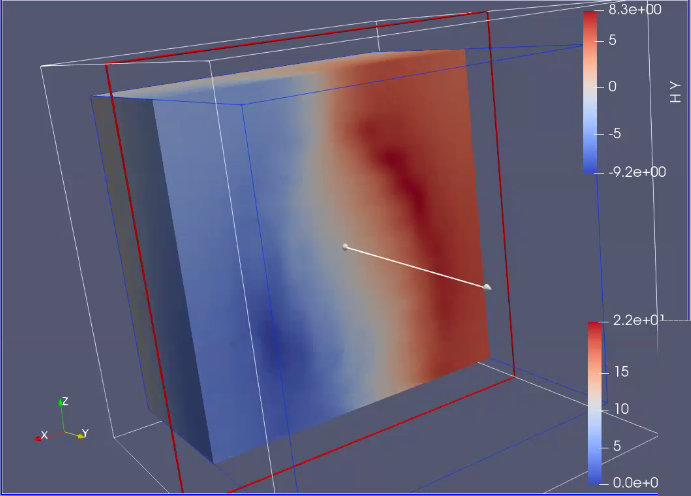}\\
	\includegraphics[scale=0.2]{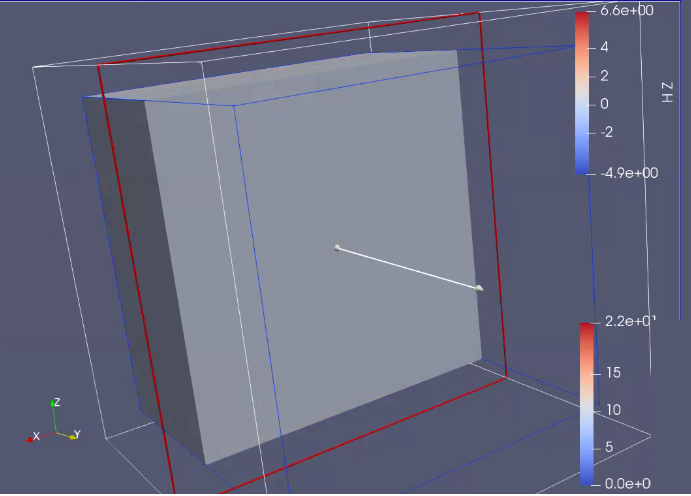}
	\includegraphics[scale=0.2]{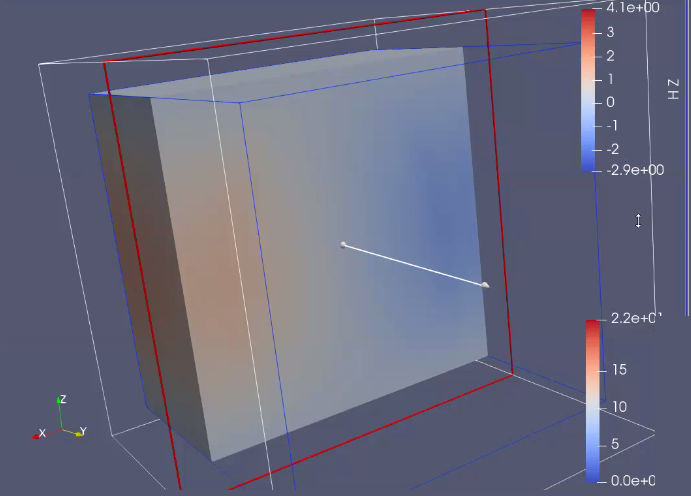}
	\includegraphics[scale=0.2]{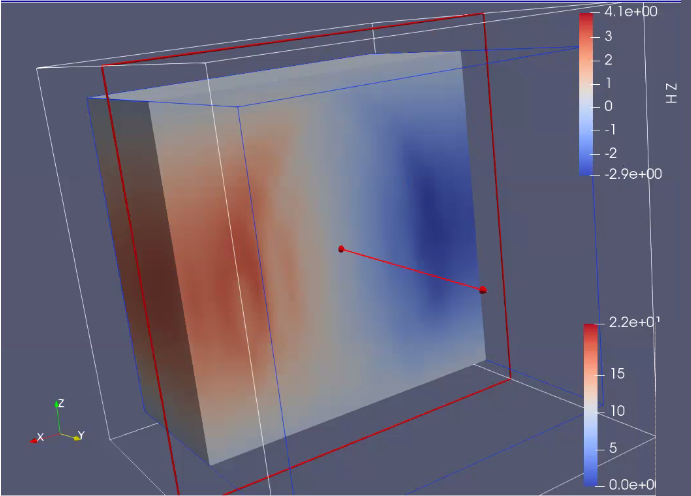}
	\includegraphics[scale=0.2]{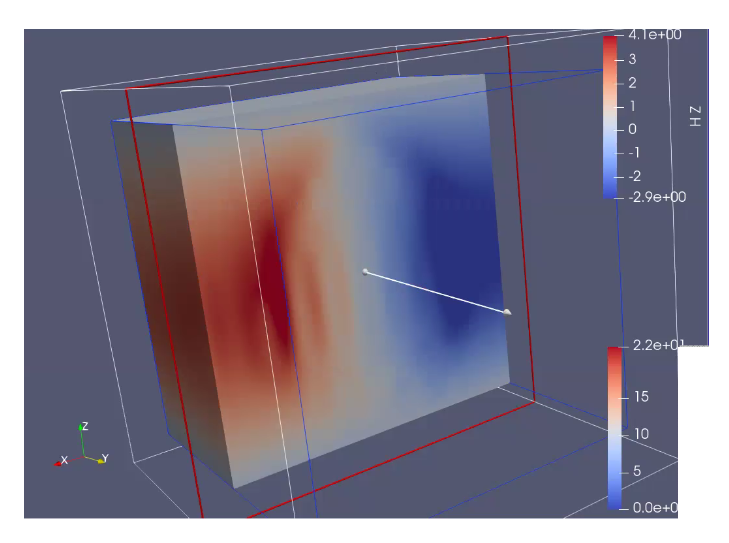}
	\includegraphics[scale=0.2]{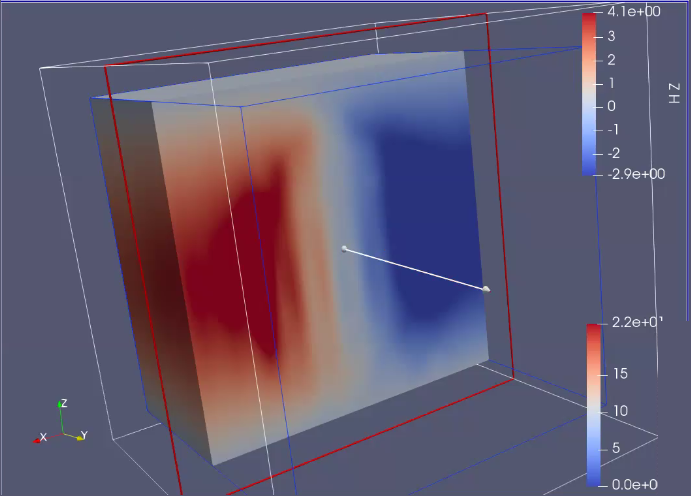}\\
	\end{center}
	\caption{Magnetic waves propagation on the human head. First row: $x$ component of the magnetic vector field. Second row: $y$ component of the magnetic vector field. Third row: $z$ component of the magnetic vector field. Columns - cross section along OYZ plane, time moments: 0, 0.25, 0.5, 0.75, 1.0s. Mesh size 32x32x32, quadratic $C^1$ B-splines.}
	\label{fig:head3}
\end{figure}

We conclude the section with numerical results presenting the electric and magnetic field over the human head data based on MRI scans.
Our simulations are based on digital data with 29 two-dimensional
slices, each one with 532 times 565 pixels. Each pixel’s intensity is a value
from the range of [0, 255], and it’s proportional to the material’s (skull, skin,
tissue, and air) normalized density. Exemplary slices of the human head from
the MRI scan are presented in Figure \ref{fig:head}.
Next, according to the MRI scan data, we used the electromagnetic waves from the manufactured solution example, with material data changing on the skull, skin, tissue, and air. 
We assume air (MRI scan data $\leq$ 1), skin or brain (tissue in general) (1 $\leq$ approximation $\leq$ 240), and skull (approximation $\geq$ 240). 
We enforce different material data using the method described in this section.
We summarize the results in Figure \ref{fig:head2} for the electric waves, and Figure \ref{fig:head3} for the magnetic waves.

For computational grids of size $8\times 8\times 8$, $16\times 16\times 16$, and $32\times 32\times 32$ with quadratic B-splines with $C^1$ continuity we check the convergence in $L^2$ and $H^1$ norms, as presented in Figure \ref{fig:head_norms}.\\

\begin{figure}
	\begin{center}
	\includegraphics[scale=0.4]{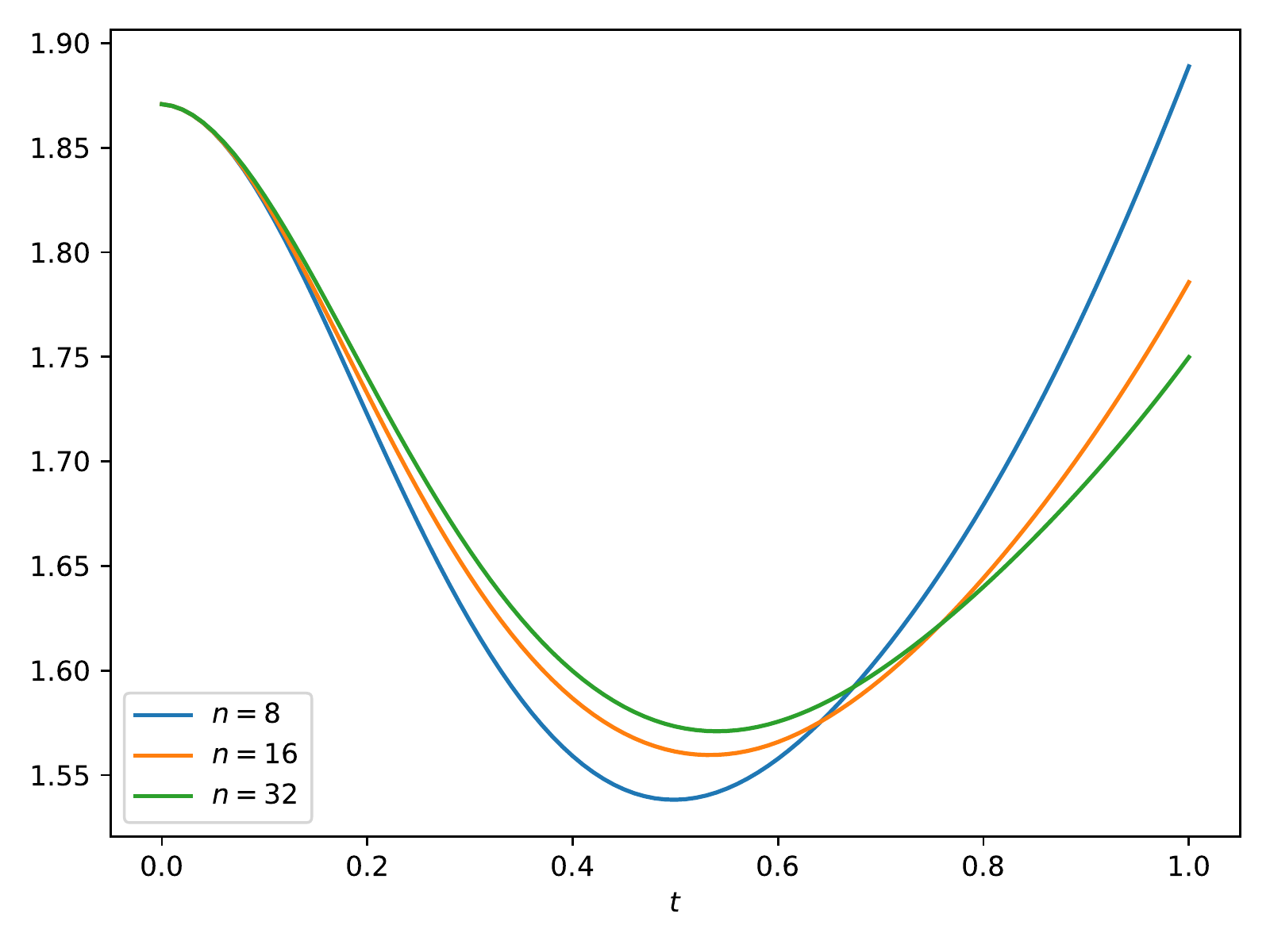}
	\includegraphics[scale=0.4]{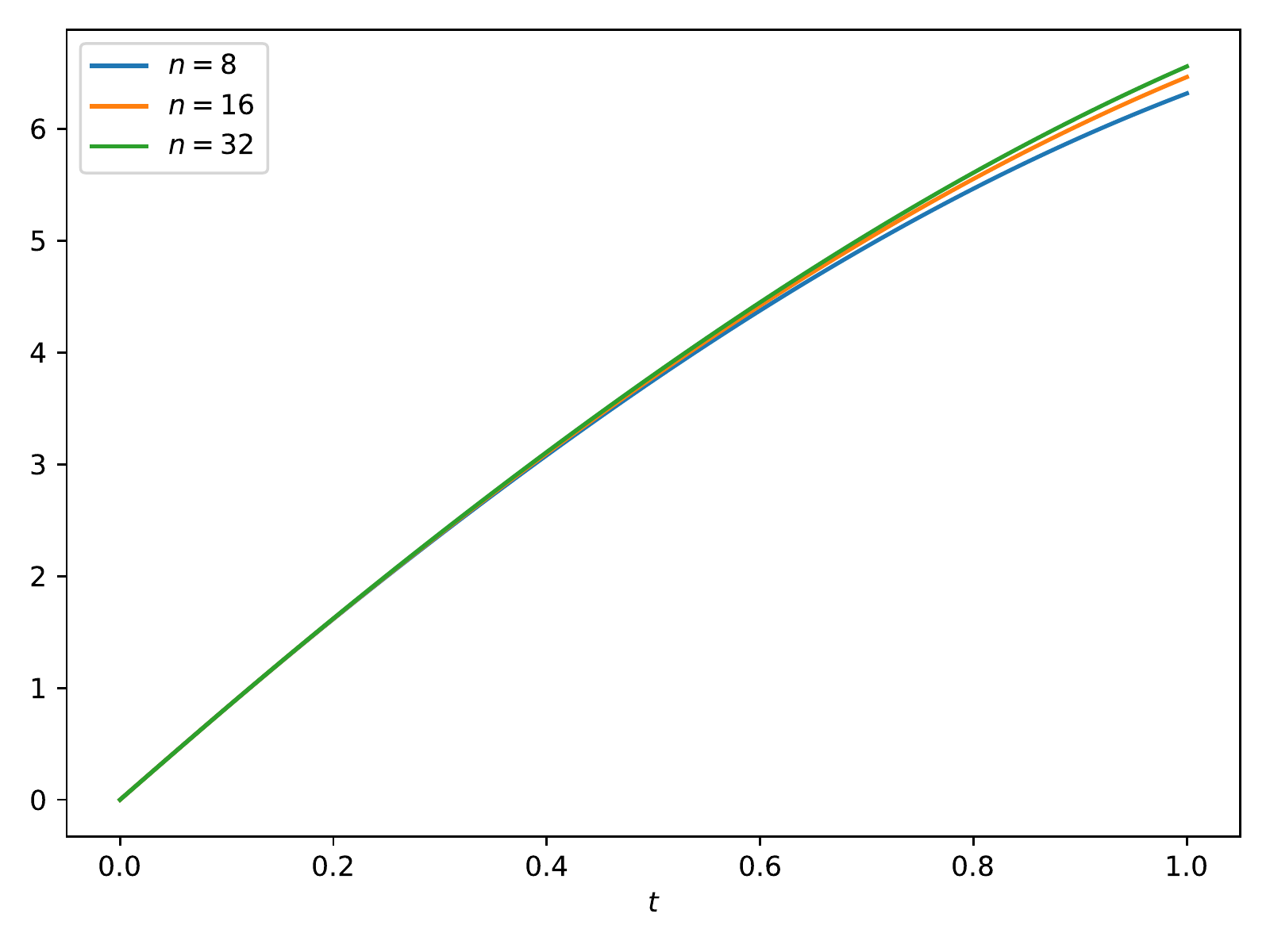}
	\end{center}
	\caption{$L^2$ norm of the electric and magnetic fields from the solutions over the human head.}
	\label{fig:head_norms}
\end{figure}

\begin{figure}
	\begin{center}
	\includegraphics[scale=0.4]{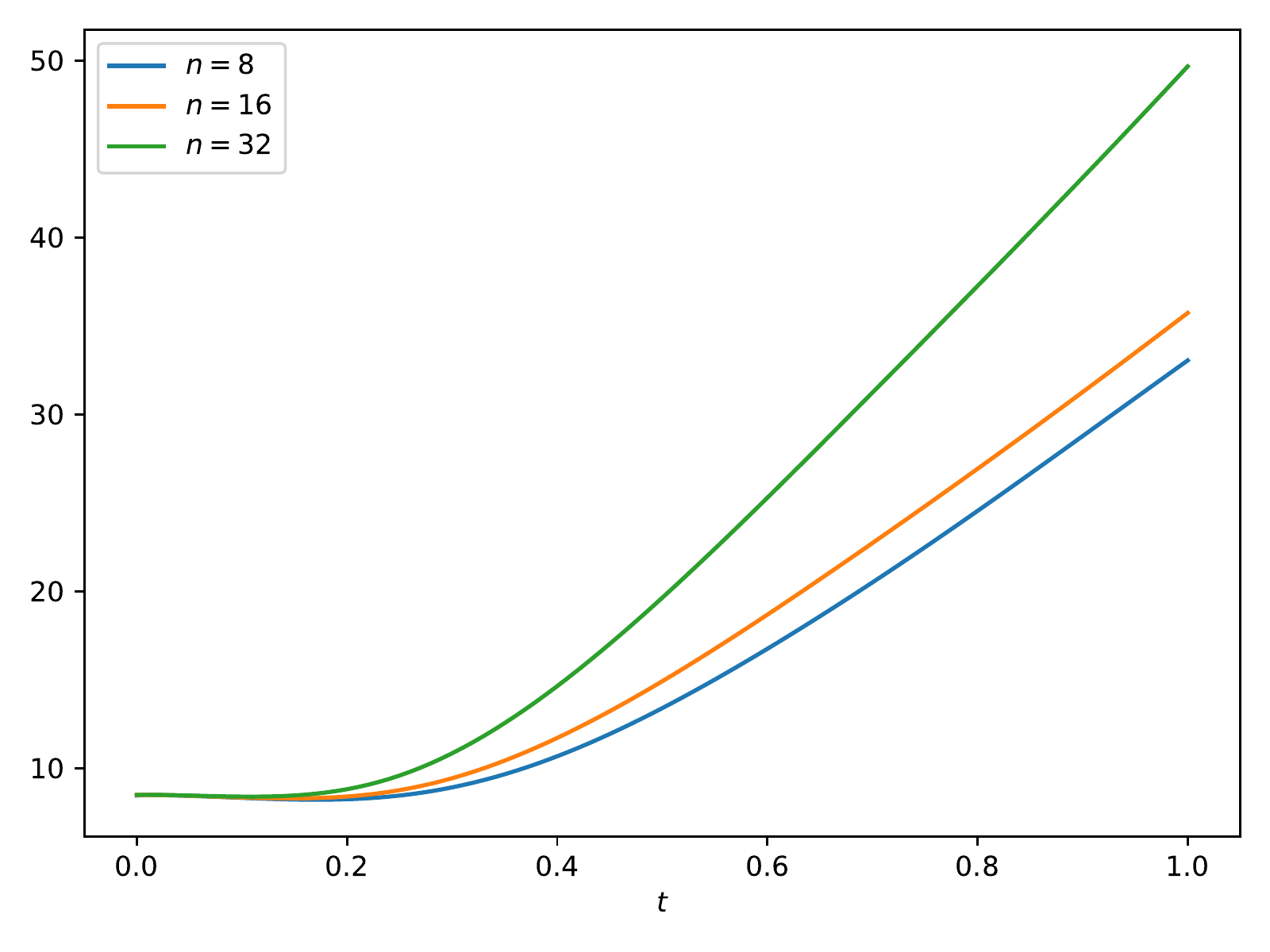}
	\includegraphics[scale=0.4]{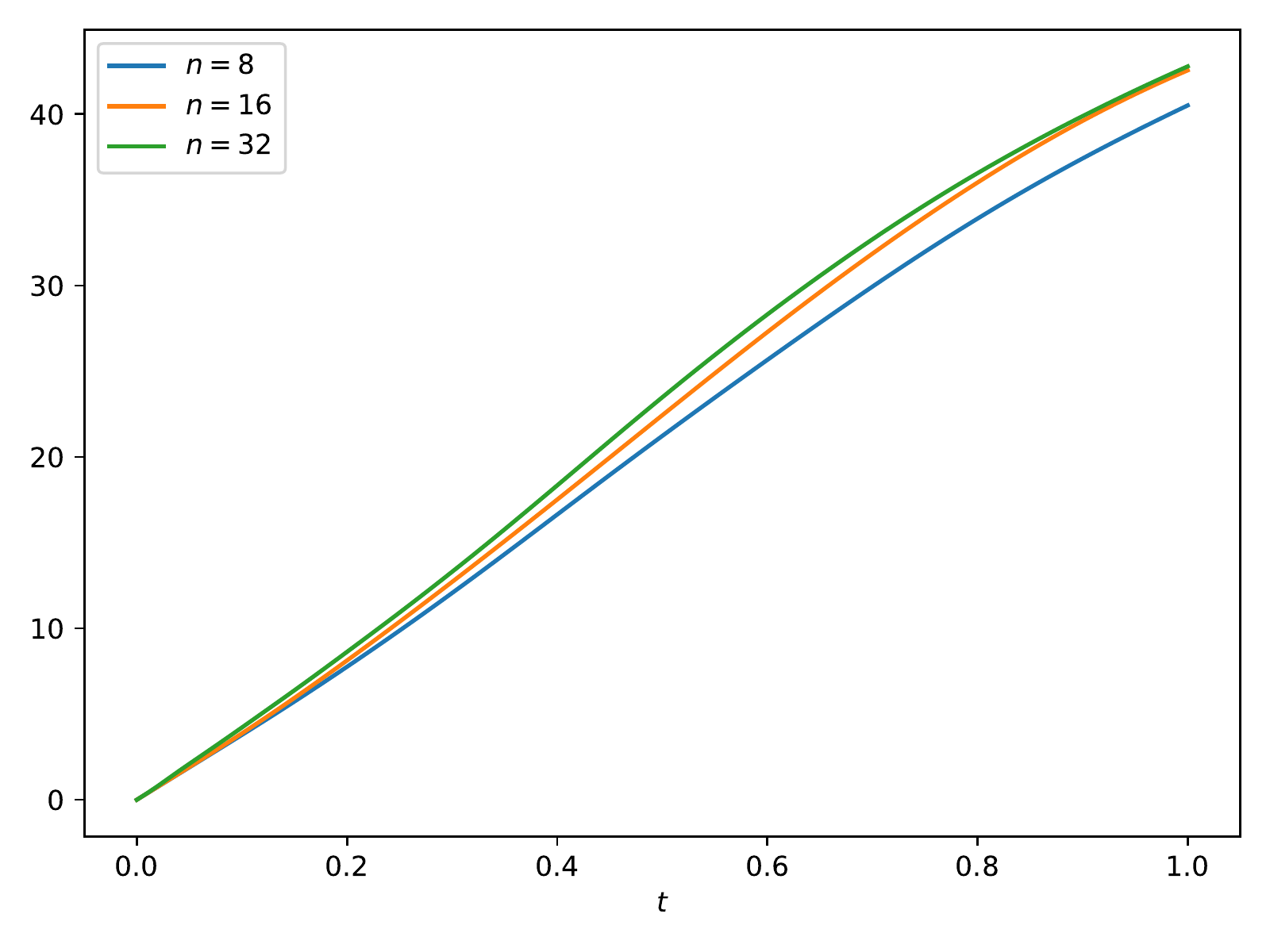}
	\end{center}
	\caption{$H^1$ norm of the electric and magnetic fields from the solutions over the human head.}
	\label{fig:head_norms}
\end{figure}

\section*{Appendix: Varying coefficients in alternating directions solver}

\begin{figure}
	\begin{center}
	\includegraphics[scale=0.4]{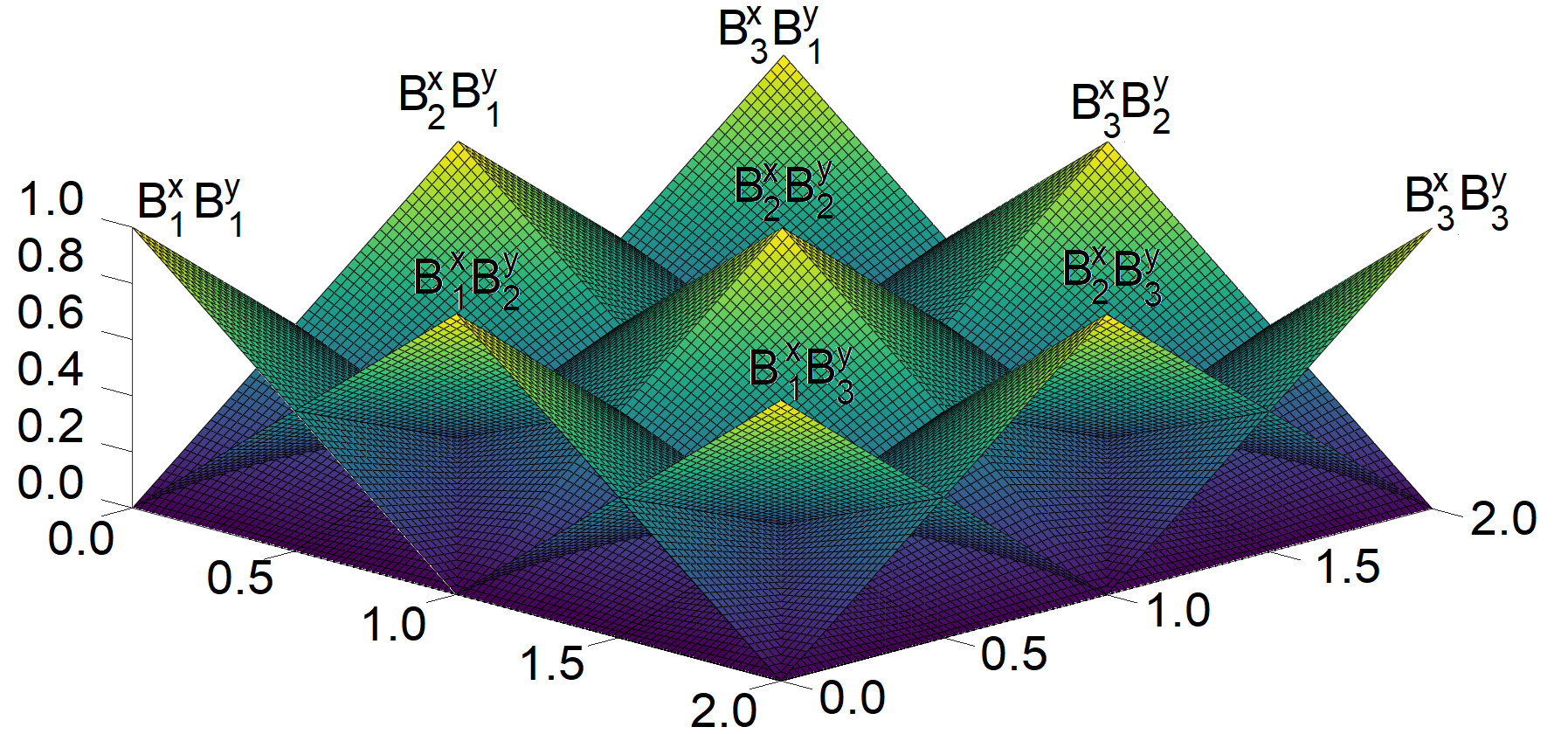}
	\end{center}
	\caption{Simple mesh with linear B-splines defined by tensor product of two knot vectors $[0 \quad 0 \quad 1\quad  2\quad  2] \times [ 0\quad  0\quad  1\quad  2\quad  2]$.}
	\label{fig:linears}
\end{figure}

This appendix shows that varying material data with test functions do not alter the linear computational cost of the direction-splitting algorithm.
To focus our attention, we consider the model projection problem augmented with  $\epsilon$ coefficients assigned to test functions (the index of the $\epsilon$ coefficients corresponds to the index of the test function in the row).
We consider a simple two-dimensional mesh with linear B-splines, presented in Figure \ref{fig:linears}. The basis is defined as a tensor product of two-knot vectors $[0 \quad 0 \quad 1\quad  2\quad  2] \times [ 0\quad  0\quad  1\quad  2\quad  2]$.

In this simple example, we employ linear B-splines. Thus some matrix entries are equal to zero (the integrals involve multiplications of B-splines that do not have common support)
\begin{eqnarray}
\begin{bmatrix}
\int \epsilon_{11} B_{1,1}(x)B_{1,1}(y)B_{1,1}(x)B_{1,1}(y)dxdy & \cdots & 
{\int \epsilon_{11} B_{3,1}(x)B_{3,1}(y)B_{1,1}(x)B_{1,1}(y)dxdy} \\
\vdots & \cdots & \vdots \\
{\int \epsilon_{33} B_{1,1}(x)B_{1,1}(y)B_{3,1}(x)B_{3,1}(y)dxdy} & \cdots & 
\int \epsilon_{33} B_{3,1}(x)B_{3,1}(y)B_{3,1}(x)B_{3,1}(y)dxdy 
\end{bmatrix} \nonumber \\
\begin{bmatrix}
u_{1,1} \\
\vdots \\
u_{3,3} 
\end{bmatrix}
= 
\begin{bmatrix}
\int F(x,y) B_{1,1}(x)B_{1,1}(y)dxdy \\
\vdots \\
\int F(x,y) B_{3,1}(x)B_{3,1}(y)dxdy 
\end{bmatrix} \nonumber \\
\end{eqnarray}

We use shorter notation $B_{1,1}(x)B_{1,1}(y)=B^x_1B^y_1$. We do not cancel out the terms that are equal to zero to illustrate the global structure of the matrix. Instead, we denote by colors the repeating terms in each of the blocks.

\begin{eqnarray}
\begin{bmatrix}
\epsilon_{11} ( \int  B^x_1\inred{B^y_1}B^x_1\inred{B^y_1} & \cdots & \int B^x_3\inred{B^y_1}B^x_1\inred{B^y_1} & \cdots & %\int B^x_1\inblue{B^y_2}B^x_1\inblue{B^y_1} & \cdots &
% \int B^x_3\inblue{B^y_2}B^x_1\inblue{B^y_1} & 
\int B^x_1\inbrown{B^y_3}B^x_1\inbrown{B^y_1} & \cdots & \int B^x_3\inbrown{B^y_3}B^x_1\inbrown{B^y_1} ) \\
\vdots & \cdots & \vdots \\
\epsilon_{31} ( \int  B^x_1\inred{B^y_1}B^x_3\inred{B^y_1} & \cdots & \int B^x_3\inred{B^y_1}B^x_3\inred{B^y_1} & \cdots & %\int B^x_1\inblue{B^y_2}B^x_3\inblue{B^y_1} & \cdots &
% \int B^x_3\inblue{B^y_2}B^x_3\inblue{B^y_1} & 
\int B^x_1\inbrown{B^y_3}B^x_3\inbrown{B^y_1} & \cdots & \int B^x_3\inbrown{B^y_3}B^x_1\inbrown{B^y_1} ) \\
\epsilon_{12} ( \int  B^x_1\incyan{B^y_1}B^x_1\incyan{B^y_2} & \cdots & \int B^x_3\incyan{B^y_1}B^x_1\incyan{B^y_2} & \cdots & %\int B^x_1\inmagenta{B^y_2}B^x_1\inmagenta{B^y_2} & \cdots &
% \int B^x_3\inmagenta{B^y_2}B^x_1\inmagenta{B^y_2} & 
\int B^x_1\ingray{B^y_3}B^x_1\ingray{B^y_2} & \cdots & \int B^x_3\ingray{B^y_3}B^x_1\ingray{B^y_2} ) \\
\vdots & \cdots & \vdots & \cdots & \vdots & \cdots & \vdots \\
\epsilon_{32} ( \int  B^x_1\incyan{B^y_1}B^x_3\incyan{B^y_2} & \cdots & \int B^x_3\incyan{B^y_1}B^x_3\incyan{B^y_2} & \cdots & %\int B^x_1\inmagenta{B^y_2}B^x_3\inmagenta{B^y_2} & \cdots &
% \int B^x_3\inmagenta{B^y_2}B^x_3\inmagenta{B^y_2} & 
\int B^x_1\ingray{B^y_3}B^x_3\ingray{B^y_2} & \cdots & \int B^x_3\ingray{B^y_3}B^x_3\ingray{B^y_2} ) \\
\epsilon_{13} ( \int  B^x_1\inred{B^y_1}B^x_1\inred{B^y_3} & \cdots & \int B^x_3\inred{B^y_1}B^x_1\inred{B^y_3} &  \cdots & % \int B^x_1\inblue{B^y_2}B^x_1\inblue{B^y_3} & \cdots &
% \int B^x_3\inblue{B^y_2}B^x_1\inblue{B^y_3} &
\int B^x_1\inbrown{B^y_3}B^x_1\inbrown{B^y_3} & \cdots & \int B^x_3\inbrown{B^y_3}B^x_1\inbrown{B^y_3} ) \\
 \vdots & \cdots & \vdots & \cdots & \vdots & \cdots & \vdots \\
\epsilon_{33} ( \int  B^x_1\inred{B^y_1}B^x_3\inred{B^y_3} & \cdots & \int B^x_3\inred{B^y_1}B^x_3\inred{B^y_3} & \cdots &  % \int B^x_1\inblue{B^y_2}B^x_3\inblue{B^y_3} & \cdots &
% \int B^x_3\inblue{B^y_2}B^x_3\inblue{B^y_3} & 
\int B^x_1\inbrown{B^y_3}B^x_3\inbrown{B^y_3} & \cdots & \int B^x_3\inbrown{B^y_3}B^x_3\inbrown{B^y_3} ) 
\end{bmatrix}
\nonumber \\
\begin{bmatrix}
u_{1,1} \\
\vdots \\
u_{3,3} 
\end{bmatrix} =
\begin{bmatrix}
\int F(x,y) B^x_1B^y_1 \\
\vdots \\
\int F(x,y) B^x_3B^y_3 
\end{bmatrix} \nonumber \\
\end{eqnarray}

Notice that each of the nine blocks (denoted by different colors) have a repeated matrix 
$$
{\cal A}=\begin{bmatrix}
\int_x  B^x_1B^x_1 & \int_x B^x_2B^x_1 & \int_x B^x_3B^x_1 \\
\int_x  B^x_1B^x_2 & \int_x B^x_2B^x_2 & \int_x B^x_3B^x_1 \\
\int_x  B^x_1B^x_3 &  \int_x B^x_2B^x_3 &\int_x B^x_3B^x_3 
\end{bmatrix}
$$

Additionally, we distinguish three different blocks, multiplied by the $\epsilon_{ij}$ constants
$$
{\cal A}_1=\begin{bmatrix}
\epsilon_{11}(\int_x  B^x_1B^x_1 & \int_x B^x_2B^x_1 & \int_x B^x_3B^x_1) \\
\epsilon_{21}(\int_x  B^x_1B^x_2 & \int_x B^x_2B^x_2 & \int_x B^x_3B^x_1) \\
\epsilon_{31}(\int_x  B^x_1B^x_3 &  \int_x B^x_2B^x_3 &\int_x B^x_3B^x_3) 
\end{bmatrix} 
$$
$$
{\cal A}_2=\begin{bmatrix}
\epsilon_{12}(\int_x  B^x_1B^x_1 & \int_x B^x_2B^x_1 & \int_x B^x_3B^x_1) \\
\epsilon_{22}(\int_x  B^x_1B^x_2 & \int_x B^x_2B^x_2 & \int_x B^x_3B^x_1) \\
\epsilon_{32}(\int_x  B^x_1B^x_3 &  \int_x B^x_2B^x_3 &\int_x B^x_3B^x_3) 
\end{bmatrix} 
$$
$$
{\cal A}_3=\begin{bmatrix}
\epsilon_{13}(\int_x  B^x_1B^x_1 & \int_x B^x_2B^x_1 & \int_x B^x_3B^x_1) \\
\epsilon_{23}(\int_x  B^x_1B^x_2 & \int_x B^x_2B^x_2 & \int_x B^x_3B^x_1) \\
\epsilon_{33}(\int_x  B^x_1B^x_3 &  \int_x B^x_2B^x_3 &\int_x B^x_3B^x_3) 
\end{bmatrix}
$$
We take them out
\begin{eqnarray}
\begin{bmatrix}
{\cal A}_{1} ( \int_y\inred{B^y_1}\inred{B^y_1}  & \int_y \inblue{B^y_2}\inblue{B^y_1} & \int_y \inbrown{B^y_3}\inbrown{B^y_1} ) \\
{\cal A}_{2} ( \int_y\incyan{B^y_1}\incyan{B^y_2} & \int_y \inmagenta{B^y_2}\inmagenta{B^y_2} & \int_y \ingray{B^y_3}\ingray{B^y_2} ) \\
{\cal A}_{3} ( \int_y \inred{B^y_1}\inred{B^y_3} &  \int_y \inblue{B^y_2}\inblue{B^y_3}  & \int_y \inbrown{B^y_3}\inbrown{B^y_3} ) 
\end{bmatrix}
\begin{bmatrix}
u_{1,1} \\
\vdots \\
u_{3,3} 
\end{bmatrix} =
\begin{bmatrix}
{\cal F}_1 \\
{\cal F}_2 \\
{\cal F}_3
\end{bmatrix} \nonumber \\
\end{eqnarray}
where we have denoted 
\begin{eqnarray}
{\cal F}_1 = \begin{bmatrix} \int F(x,y) B^x_1B^y_1  \\ \int F(x,y) B^x_2B^y_1 \\ \int F(x,y) B^x_3B^y_1 \end{bmatrix} \quad
{\cal F}_2 = \begin{bmatrix} \int F(x,y) B^x_1B^y_2  \\ \int F(x,y) B^x_2B^y_2 \\ \int F(x,y) B^x_3B^y_2 \end{bmatrix} \quad
{\cal F}_3 = \begin{bmatrix} \int F(x,y) B^x_1B^y_3  \\ \int F(x,y) B^x_2B^y_3 \\ \int F(x,y) B^x_3B^y_3 \end{bmatrix} 
\end{eqnarray}
We multiply blocks by ${\cal A}_1^{-1}$,  ${\cal A}_2^{-1}$, and ${\cal A}_3^{-1}$, and we define 
\begin{eqnarray}
{\cal G}_1={\cal A}_1^{-1} {\cal F}_1 \Longleftrightarrow {\cal A}_1{\cal G}_1= {\cal F}_1 \nonumber \\
{\cal G}_2={\cal A}_2^{-1} {\cal F}_2 \Longleftrightarrow {\cal A}_2{\cal G}_2= {\cal F}_2 \nonumber \\
{\cal G}_3={\cal A}_3^{-1} {\cal F}_3\Longleftrightarrow {\cal A}_3{\cal G}_3={\cal F}_3 \label{eq:split1}
\end{eqnarray}
Finally, we have
\begin{eqnarray}
\begin{bmatrix}
 \int_y\inred{B^y_1}\inred{B^y_1}  & \int_y \inblue{B^y_2}\inblue{B^y_1} & \int_y \inbrown{B^y_3}\inbrown{B^y_1}  \\
 \int_y\incyan{B^y_1}\incyan{B^y_2} & \int_y \inmagenta{B^y_2}\inmagenta{B^y_2} & \int_y \ingray{B^y_3}\ingray{B^y_2}  \\
 \int_y \inred{B^y_1}\inred{B^y_3} &  \int_y \inblue{B^y_2}\inblue{B^y_3}  & \int_y \inbrown{B^y_3}\inbrown{B^y_3} 
\end{bmatrix}
\begin{bmatrix}
u_{1,1} & u_{1,2} & u_{1,3} \\
u_{2,1} & u_{2,2} & u_{2,3} \\
u_{3,1} & u_{3,2} & u_{3,3} \\
\end{bmatrix}
 =
\begin{bmatrix}
{\cal G}_1 & 
{\cal G}_2 &
{\cal G}_3
\end{bmatrix}  \label{eq:split2}
\end{eqnarray}

Both systems (\ref{eq:split1})-(\ref{eq:split2}) can be solved in a linear computational cost due to the banded structures of mass matrices build with one-dimensional B-splines.

\section*{Conclusions}
\label{sec:conc}

In this paper, we applied the isogeometric analysis (IGA) to discretize the time-dependent Maxwell equations. Furthermore, we used the alternating direction splitting with an implicit time integration scheme for fast solution (ADI solvers).
Our method delivers a linear computational cost ${\cal O}(N)$ solver, unconditional stability in time, delivered by the implicit time integration scheme, and the second-order accurate time integration scheme. Additionally, we showed that the linear cost of the solver is preserved, even if we vary material data over the computational domain.
We mix benefits of the state-of-the-art modern methods, namely 
the Isogeometric Finite Element Method (IGA-FEM) \cite{p1}, and Alternating Direction Implicit solvers (ADI) \cite{ADS1,ADS2,maxwell1}.
We showed how to run our linear computational cost solver on non-regular material data, such as the human head's tissue and skull. Our method allows for fast and reliable simulations of electromagnetic waves propagation through non-regular biological tissues.

\section*{Acknowledgments}
This work is supported by National Science Centre, Poland grant no. 2017/26/M/ ST1/ 00281.

% Either type in your references using
% \begin{thebibliography}{}
% \bibitem{}
% Text
% \end{thebibliography}
%
% or
%
% Compile your BiBTeX database using our plos2015.bst
% style file and paste the contents of your .bbl file
% here. See http://journals.plos.org/plosone/s/latex for 
% step-by-step instructions.
% 

\end{document}